\documentclass[bj]{imsart}

\RequirePackage{amsthm,amsmath,amsfonts,amssymb}
\RequirePackage[numbers]{natbib}
\RequirePackage[colorlinks,citecolor=blue,urlcolor=blue]{hyperref}
\RequirePackage{graphicx}

\startlocaldefs
\numberwithin{equation}{section}

\newtheorem{theorem}{Theorem}[section]
\newtheorem{lemma}[theorem]{Lemma}
\theoremstyle{remark}

\newtheorem{assumption}[theorem]{Assumption}
\newtheorem{remark}[theorem]{Remark}

 \newcommand{\dif}{\mathrm{d}}
 
 \newcommand{\A}{\mathcal{A}}

 \newcommand{\Y}{\mathcal{Y}}

 \newcommand{\W}{\mathcal{W}}           
 \newcommand{\E}{\mathbb{E}}            
 
 \newcommand{\e}{\varepsilon}
 
 \newcommand{\p}{\partial}

 \newcommand{\Ll}{\langle}
 \newcommand{\Rr}{\rangle}

 \newcommand{\R}{\mathbb{R}}
 
 \newcommand{\PP}{\mathbb{P}}
 \newcommand{\mcl}{\mathcal}
 \newcommand{\vp}{\varphi}
 \newcommand{\Be}{\begin{equation}}
 \newcommand{\Ee}{\end{equation}}
  \newcommand{\Bes}{\begin{equation*}}
 \newcommand{\Ees}{\end{equation*}}
  \newcommand{\Bey}{\begin{eqnarray}}
 \newcommand{\Eey}{\end{eqnarray}}
 \newcommand{\Beys}{\begin{eqnarray*}}
 \newcommand{\Eeys}{\end{eqnarray*}}
 \newcommand{\BT}{\begin{thm}}
 \newcommand{\ET}{\end{thm}}
 \newcommand{\Bp}{\begin{proof}}
 \newcommand{\Ep}{\end{proof}}
 \newcommand{\BL}{\begin{lem}}
 \newcommand{\EL}{\end{lem}}
 \newcommand{\BP}{\begin{proposition}}
 \newcommand{\EP}{\end{proposition}}
 \newcommand{\BC}{\begin{corollary}}
 \newcommand{\EC}{\end{corollary}}
 \newcommand{\BR}{\begin{rem}}
 \newcommand{\ER}{\end{rem}}
 \newcommand{\BD}{\begin{defn}}
 \newcommand{\ED}{\end{defn}}
 \newcommand{\BI}{\begin{itemize}}
 \newcommand{\EI}{\end{itemize}}
 
 \newcommand{\tl}{\tilde}

\endlocaldefs

\begin{document}

\begin{frontmatter}
\title{Central limit theorem and Self-normalized Cram\'er-type moderate deviation for Euler-Maruyama Scheme}

\begin{aug}
\author[A,B]{\fnms{Jianya} \snm{Lu}\ead[label=e1]{jianya.lu@connect.um.edu.mo}},
\author[C]{\fnms{Yuzhen} \snm{Tan}\ead[label=e3]{tanyuzhensdu@gmail.com}}
\and
\author[A,B]{\fnms{Lihu} \snm{Xu}\ead[label=e2,mark]{lihuxu@umac.mo}}
\address[A]{Department of Mathematics, Faculty of Science and Technology, University of Macau, Macau, China.
\printead{e1,e2}}

\address[B]{UM Zhuhai Research Institute, Zhuhai, China.}

\address[C]{Zhongtai Securities Institute for Financial Studies, Shandong University, Jinan, China
\printead{e3}}

\end{aug}

\begin{abstract}
We consider a stochastic differential equation and its Euler-Maruyama (EM) scheme, under some appropriate conditions, they both admit a unique invariant measure, denoted by $\pi$ and $\pi_\eta$ respectively ($\eta$ is the step size of the EM scheme). We construct an empirical measure $\Pi_\eta$ of the EM scheme as a statistic of $\pi_\eta$, and use Stein's method developed in \citet{FSX19} to prove a central limit theorem of $\Pi_\eta$. The proof of the self-normalized Cram\'er-type moderate deviation (SNCMD) is based on a standard decomposition on Markov chain, splitting $\eta^{-1/2}(\Pi_\eta(.)-\pi(.))$ into a martingale difference series sum $\mcl H_\eta$ and a negligible remainder $\mcl R_\eta$. We handle $\mcl H_\eta$  by the time-change technique for martingale, while prove that $\mcl R_\eta$ is exponentially negligible by concentration inequalities, which have their independent interest. Moreover, we show that SNCMD holds for  $x = o(\eta^{-1/6})$, which has the same order as that of the classical result in \citet{shao1999cramer,JSW03}.
\end{abstract}


\end{frontmatter}

\section{Introduction}
We consider the following stochastic differential equation (SDE) on $\R^d$:
\begin{align}\label{SDE}
\dif X_{t}=g(X_{t}) \dif t+\sigma(X_t) \dif B_{t},\quad X_0=x,
\end{align}
where $\sigma: \R^d \rightarrow \R^{d \times d}$ and $g: \R^d \rightarrow \R^d$ satisfy {\bf Assumption} \ref{assu} below,  and $B_t$ is a $d-$dimentional standard Brownian motion. Given a step size $\eta$, the Euler-Maruyama (EM) scheme of \eqref{SDE} reads as
\begin{align}\label{Lang0}
\theta_{k+1}=\theta_k+\eta g(\theta_k)+\sqrt{\eta}\sigma(\theta_k) \xi_{k+1}, \ \ \ \ k \ge 0,
\end{align}
where $(\xi_k)_{k\ge1}$ are i.i.d. standard $d$-dimensional normal random vectors. When $g$ and $\sigma$ are both Lipschitz, \eqref{SDE} admits a unique strong solution and the following strong approximation error bound holds, see \citet{Mao97} : for any $T>0$,
\Be  \label{e:SAppr}
\E |X_{T}-\theta_{[T/\eta]}|^2 \le C_T \eta,
\Ee
the constant $C_T$ usually tends to $\infty$ as $T \to \infty$ and $[x]$ denotes the integer part of $x$ for a $x>0$. When $g$ or $\sigma$ is irregular, there have recently been some works, see \citet{Bao19} for the convergence rate of degenerate SDEs. We refer the reader to \citet{Bao18,Shao18} for the EM scheme of path-dependent SDEs and to \citet{Bao12} for that of stochastic differential delay equations.

Let us first discuss a special case of \eqref{SDE} in which $\sigma(x) \equiv I_{d \times d}$, $d\times d$ identity matrix, and $g(x)=-\nabla U(x)$ with $U$ being a potential, it is well known that \eqref{SDE} is a gradient system and admits a unique ergodic measure $\pi$ proportional to $e^{-U(x)}$ from \citet{Tweedie1996}. \eqref{Lang0} is called unadjusted Langevin algorithm (ULA) with constant step size, \citet{Tweedie1996} mainly established some criteria for the ergodicity of $\theta_k$, while \citet[Theorem 2]{Dal17} gave an explicit error in total variation distance between $\theta_k$ and $\pi$ in terms of $d,k,\eta$ when $\nabla U$ is Lipschitz and strong convex.
Replacing the strong convexity assumption in \citet{Dal17} with a strong convexity at infinity condition, \citet{Majka2020} used a coupling method to show the Wasserstein-2 distance between $\theta_k$ and $\pi$  were bounded by $C[(1-\eta)^{k/2}+\eta^{\frac14}]$. When $\nabla U$ is third order differentiable with a appropriate growth condition but not necessarily Lipschitz, \citet{FSX19} showed that as long as the above \eqref{Lang0} admits a unique ergodic measure $\pi_\eta$, then the Wasserstein-1 distance between $\pi_\eta$ and $\pi$ is bounded by $\sqrt{\eta}$ up to a logarithmic correction.
 For more research about Langevin algorithm, we refer the reader to \citet{JY1,JY2,JY3} and the references therein.

The motivations of studying the central limit theorem (CLT) and the self-normalized Cram\'er type moderate deviation (SNCMD) of ULA are two folds. One is that there have been many central limit theorems and moderate deviation results for Markov chain Monte Carlo (MCMC) algorithm, see \citet{DuJo17,MeTw93b,Moral15,Nyq17,Ti94}, whereas there are very few these type of fluctuation theorems for Langevin algorithm. The other is that our result provides a new example for SNCMD for dependent time series, and also a new example that applies Stein's method to prove SNCMD, see \citet{CFS13,SZZ20+}. Note that there are not many results for SNCMD for dependent time series, see \citet{CSWX16, fan2020cramer, Fan2019Self, Fan2020, FLS20,jing2015cramer,shao2016cramer} and the references therein.

Let us briefly describe our main results and methods as follows. We construct an empirical measure $\Pi_\eta$ as a statistic of the ergodic measure $\pi_\eta$ of \eqref{Lang0}, for any function $h \in C^2_b(\R^d,\R)$ (see the definition of $C^2_b(\R^d,\R)$ below), we study the CLT and SNCMD of $\Pi_\eta(h)$. In order to prove the CLT, we apply Stein's method developed in \citet{FSX19}. {\bf Assumption} \ref{assu} guarantees that \eqref{Lang0} admits a unique invariant measure $\pi_\eta$, while the restriction of $h \in C^2_b(\R^d,\R)$ ensures that the solution $\varphi$ of Stein's equation \eqref{stein} has bounded 4th order derivatives. Note that the ergodicity of \eqref{SDE} does not imply that of \eqref{Lang0}, see \citet{Tweedie1996}. The proof of SNCMD is based on a standard decomposition on Markov chain, splitting $\eta^{-1/2}(\Pi_\eta(h)-\pi(h))$ into a martingale difference series sum $\mcl H_\eta$ and a negligible remainder $\mcl R_\eta$. We handle $\mcl H_\eta$  by the time-change technique for martingale, while prove that $\mcl R_\eta$ is exponentially negligible by concentration inequalities, which have their independent interest. Moreover, we show that SNCMD holds for $x = o(\eta^{-1/6})$, which has the same order as that of the classical result in \citet{shao1999cramer,JSW03}. Indeed, the limit $\lim_{\eta \rightarrow 0} (\Pi_\eta(h)-\pi(h))=0$ can be understood as a law of large number (LLN), after zooming in on it by a scale $\eta^{-1/2}$, $\eta^{-1/2}(\Pi_\eta(h)-\pi(h))$ has a normal distributed fluctuation. Our result showed that this fluctuation is uniformly comparable with normal distribution for all $x \in (c\eta^{1/6},o(\eta^{-1/6}))$. In contrast, Shao et al.'s result means that by zooming in on $\frac 1n \sum_{i=1}^n X_i-\E X$ with a scale $n^{1/2}$,
$n^{1/2}(\frac 1n \sum_{i=1}^n X_i-\E X)$ has a normal distributed fluctuation uniformly comparable with normal distribution for all $x \in [0,o(n^{1/6}))$.

The paper is organized as the following. Our main results are stated and discussed in Section 2. In Section 3, we provide some preliminary lemmas. The proof of the CLT is given in Section 4. In Section 5, we give the proof of SNCMD. The details of the proof of preliminary lemmas are deferred to Appendix.

We finish this section by introducing some notations which will be frequently used in sequel. For $x\in\R^d$, $x_i$ denotes the $i-$th element of $x$. For function $f:\R^d\to\R$, denote $\nabla_{i,j,k}^3f(x)=\frac{\partial^3 f}{\partial x_i\partial x_j\partial x_k}$ with $i,j,k=1,2,...,d$. $C_b^k(\R^d,\R)$ with $k \ge 1$ denotes the collection of all bounded $k-$th order continuously differentiable
functions. The symbols $C$ and $c$ denote positive numbers depending on $g$ and $\sigma$, $C_p$ and $c_p$ denote positive numbers depending on $g$, $\sigma$ and the parameter $p$. Their values may vary from line to line.  We denote the Euclidean norm of $\R^d$ by $|\cdot|$ and for higher rank tensors by $\|\cdot\|$. For function $f$, we denote $\|f\|=\sup_{x\in\R^d}\|f(x)\|$. If a random variable $\xi$ has a probability distribution $\mu$, we write $\xi \sim \mu$. Let $\{a_n\}_{n \ge 1}$ and $\{b_n\}_{n \ge 1}$ be two nonnegative real number sequences, if there exist some $C>0$ such that $|a_n| \le C b_n$, we write $a_n=O(b_n)$. If $\lim_{n \rightarrow \infty} \frac{a_n}{b_n}=0$, we write $a_n=o(b_n)$.

\section{Main results}\label{s:2}
\begin{assumption}\label{assu}
$\sigma(x) \equiv \sigma$ with $\sigma$ being an invertible $d \times d$ matrix.
$g: \R^{d} \rightarrow \R^{d}$ is second order differentiable.  There exist $L,K_1>0$ and $K_2\ge0$ such that for every $x, y\in\mathbb{R}^d$
    \begin{align}\label{assu1}
	|g(x)-g(y)|\leq L|x-y|,
	\end{align}
	\begin{align}\label{assu2}
	\Ll g(x)-g(y),x-y\Rr\le -K_1|x-y|^2+K_2.
	\end{align}
Moreover, the second order derivative of $g$ is bounded.
\end{assumption}
\begin{remark}
	It is easy to see that the assumption (\ref{assu1}) implies
	\begin{align}\label{assu3}
	|g(x)|^2\le 2L^2|x|^2+2|g(0)|^2, \ \ \ \|\nabla g\| \le L,
	\end{align}
	and that the assumption (\ref{assu2}) and Young's inequality imply
	\begin{eqnarray}\label{assu4}
	\Ll x,g(x)\Rr&=& \Ll x-0,g(x)-g(0)\Rr+\Ll x,g(0)\Rr\nonumber\\
    &\le&-K_1|x|^2+K_2+\frac{K_1}{2}|x|^2+\frac{1}{2K_1}|g(0)|^2
	=-\frac{K_1}{2}|x|^2+C.
	\end{eqnarray}
The condition that $g$ has bounded second order derivative is only needed for proving the regularity to the solution of Stein's equation.
\end{remark}
Under {\bf Assumption} \ref{assu}, the Euler-Maruyama scheme reads as
\begin{align}\label{Lang}
\theta_{k+1}=\theta_k+\eta g(\theta_k)+\sqrt{\eta}\sigma \xi_{k+1},\ \ \ k \ge 0,
\end{align}
where $\theta_0=x$ and $(\xi_k)_{k\ge1}$ are i.i.d. standard $d$-dimensional normal random vectors.
\begin{lemma}  \label{l:Erg2}
	Under Assumption \ref{assu}, SDE (\ref{SDE}) and $(\theta_k)_{k\ge0}$ are both ergodic with invariant measures $\pi$ and $\pi_\eta$ respectively.
\end{lemma}
\begin{proof}
	The proof will be given in Appendix A.	
\end{proof}

The generator $\mathcal{A}$ of (\ref{SDE}) is given by
\begin{align}\label{generator}
\mathcal{A}f(x)=\Ll g(x),\nabla f(x)\Rr+\frac12\Ll \sigma\sigma^{\mathrm{T}}, \nabla^2 f(x)\Rr_{\mathrm{HS}},
\end{align}
where $\mathrm{T}$ is the transpose operator and $\Ll A,B\Rr_{\mathrm{HS}}:=\sum_{i,j=1}^dA_{ij}B_{ij}$ for $A,B\in\R^{d\times d}$, and $f \in C^2_b(\R^d,\R)$.
To approximate the behavior of $(X_t)_{t\ge0}$, we can use the Euler-Maruyama scheme to discrete $(\ref{SDE})$.

For a small $\eta \in (0,1)$, define
 \begin{eqnarray}  \label{e:PiSta}
	\Pi_\eta(\cdot)&=&\frac{1}{[\eta^{-2}]}\sum_{k=0}^{[\eta^{-2}]-1} \delta_{\theta_k}(\cdot),
\end{eqnarray}
where $\delta_y(\cdot)$ is a delta measure of $y$, i.e., for any $A \subset \R^d$, $\delta_y(A)=1$ if $y \in A$ and $\delta_y(A)=0$ if $y \notin A$. We shall see that $\Pi_\eta$ is an asymptotically consistent statistic of $\pi$ as $\eta \rightarrow 0$.

Parallel to the CLT and tail probability estimates of MCMC algorithms, see \citet{Rob04}, it is natural to consider those for
$\Pi_\eta$. For a test function $h:\R^d \rightarrow \R$, we consider the limit of
$\frac{\Pi_\eta(h)-\pi(h)}{\sqrt{\eta}}$  with $\pi(h)=\int_{\R^d} h(x)\pi(\dif x)$.
Our first main result is
\begin{theorem}\label{thm-CLT}
	Suppose that Assumption \ref{assu} holds. Let $h \in C^2_b(\R^d,\R)$, then we have
	\begin{align}
	\frac 1{\sqrt{\eta}}(\Pi_\eta(h)-\pi(h))\Rightarrow N(0,\pi(|\sigma^{\rm T}\nabla\vp|^2)),\quad\text{as }\eta\to0,
	\end{align}
where $\vp$ is the solution to the following Stein's equation:
	\begin{align}\label{stein}
	h-\pi(h)=\A \vp,
	\end{align}
and $\mathcal{A}$ is the generator (\ref{generator}) of the SDE \eqref{SDE}.
\end{theorem}

Let $\E_k[\cdot]$ and $\mathbb{P}_k(\cdot)$ be respectively the conditional expectation $\E[\cdot|\theta_k]$ and conditional probability $\mathbb{P}(\cdot|\theta_k)$.  Let $\Phi(x)$ be the standard normal distribution function. Denote  \footnote{Prof. Fuqing Gao suggested that we replace the self-normalized factor $\frac{1}{[\eta^{-2}]}\sum_{k=0}^{[\eta^{-2}]-1}\langle\sigma^{\mathrm{T}}\nabla \vp(\theta_k ),\xi_{k+1}\rangle^2$ in the previous version by $\frac{1}{[\eta^{-2}]}\sum_{k=0}^{[\eta^{-2}]-1}|\sigma^{\mathrm{T}}\nabla \vp(\theta_k )|^2$. Since $(\theta_k)_{k \ge 0}$ is observable whereas $\{\xi_k\}_{k \ge 1}$ is not known, the new self-normalized factor is more natural. }
\begin{align*}
\Y_\eta =\frac{1}{[\eta^{-2}]}\sum_{k=0}^{[\eta^{-2}]-1}|\sigma^{\mathrm{T}}\nabla \vp(\theta_k )|^2,
\quad
\W_\eta= \frac{\eta^{-\frac 12}(\Pi_\eta(h)- \pi(h))}{\sqrt{{\Y_\eta }}}.
\end{align*}

Our second main result is the SNCMD of $\W_\eta$ as follows.

\begin{theorem}\label{thm-self-normal}
Suppose that Assumption \ref{assu} holds. Let $\theta_0\sim\pi_\eta$ and $h \in C^2_b(\R^d,\R)$, we have
\begin{equation}\label{eq-self-normal}
\frac{\PP(\W_\eta \geq x)}{1-\Phi(x)} = 1+O\big( x\eta^{\frac16}+\eta^{\frac16}\big)
\end{equation}
uniformly for $c\eta^{\frac16}\le x=o\big(\eta^{-\frac 16}\big)$ as $\eta$ vanishes, where $c$, $O$ and $o$ depend on $L,K_1,K_2$,$|g(0)|^2$, $\sigma$.
\end{theorem}

For the simplicity of notations below, without loss of generality, we assume from now on that $\eta \in (0,1)$ is a small number such that $\eta^{-1}$ is an integer. We also denote
$$m=\eta^{-2}$$ and often write $\eta^{-1}$ as $m\eta$ for notational simplicity. Denote
$$\Delta \theta_k=\theta_{k+1}-\theta_k, \ \ \ \ k \ge 0.$$

\section{Auxiliary Lemmas for Theorem \ref{thm-CLT} and Theorem \ref{thm-self-normal}}
\subsection{The strategy of proving Theorem \ref{thm-CLT} and Theorem \ref{thm-self-normal}}
The strategy of proving Theorem \ref{thm-CLT} and Theorem \ref{thm-self-normal} is to decompose $\eta^{-\frac 12}(\Pi_\eta(h)-\pi(h))$ into a martingale and a remainder as in \eqref{e:Decompose} below, showing that the remainder is negligible, while the martingale converges weakly to a normal distribution and satisfies the SNCMD. This type of decomposition is typical for proving CLT for semi-martingales, see e.g., \citet[(28)]{Teh2016}.
\begin{lemma}\label{solpoisson}
Let $h\in C_b^2(\R^d,\R)$, a solution to Stein's equation (\ref{stein}) is given by
\begin{align*}
  \vp(x)=-\int_0^\infty\E[h(X_t(x))-\pi(h)]\dif t,
\end{align*}
where $X_t(x)$ is the solution of equation (\ref{SDE}) with initial value $x$. Moreover,
\begin{align*}
  \|\nabla^k \vp\|\le C,\quad k=0, 1,2,3,4.
\end{align*}
\end{lemma}
\begin{proof}
Denote $\hat h=h-\pi(h)$ and $P_th(x)=\E[h(X_t(x))]$. Following the exponential ergodicity of $\{X_t\}_{t\ge0}$, i.e. (\ref{ergodic1}), one has
$$|\int_0^\infty P_s \hat h(x)\dif s|\le \int_0^\infty|P_s \hat h(x)|\dif s\le CV(x)\int_0^\infty e^{-cs}\dif s<\infty.
$$
Thus $\int_0^\infty P_s \hat h(x)\dif s$ is well defined. For any $\e>0$, it is known that $\e-\A$ is invertible (cf. \citet[pp. 158-159]{applebaum2009levy}), and
$$(\e-\A)^{-1}\hat h=\int_0^\infty e^{-\e t}P_t\hat h \dif t,$$
i.e.,
$$\e\int_0^\infty e^{-\e t}P_t\hat h \dif t-\hat h=\A\left(\int_0^\infty e^{-\e t}P_t\hat h \dif t\right).$$
Let $\e\to0+$,
$$\e\int_0^\infty e^{-\e t}P_t\hat h \dif t-\hat h\to-\hat h,\quad \int_0^\infty e^{-\e t}P_t\hat h \dif t\to\int_0^\infty P_t\hat h \dif t.$$
Since $\A$ is a closed operator (cf. \citet[Theorem 2.2.6]{partington04}), $\int_0^\infty P_t\hat h \dif t$ is in the domain of $\A$ and
$$\hat h=\mcl A\left(-\int_0^\infty P_t\hat h \dif t\right).$$

By \citet[Theorem 2.6]{KP08}, we know that $\varphi \in C^3_b(\R^d,\R)$. Denoting $\varphi_i=\p_{x_i} \varphi$ for $i=1,...,d$, it satisfies
$$\mcl A \varphi_i=\p_{x_i} h-\p_{x_i} g \varphi,$$
it is easy to check that the right hand side of this equation belongs to $C^1_b(\R^d)$ by Assumption \ref{assu}, we know that $\varphi_i \in C^3_b(\R^d,\R)$ by \citet[Theorem 2.6]{KP08}. Hence, $\varphi \in C_b^4(\R^d,\R)$.
\end{proof}

By Stein's equation (\ref{stein}), we have,
\begin{align*}
\Pi_\eta(h)-\pi(h)&=\frac{1}{m}\sum_{k=0}^{m-1}\left(h(\theta_k )-\pi(h)\right)=\frac{1}{m}\sum_{k=0}^{m-1}\mathcal{A}\vp(\theta_k )\\
&=\eta\sum_{k=0}^{m-1}\left[\mathcal{A}\vp(\theta_k )\eta-\left(\vp(\theta_{k+1} )-\vp(\theta_k )\right)\right]+\eta\sum_{k=0}^{m-1}\left(\vp(\theta_{k+1} )-\vp(\theta_k )\right)\\
&=\eta[\vp(\theta_m )-\vp(\theta_0 )]+\eta\sum_{k=0}^{m-1}\left[\mathcal{A}\vp(\theta_k )\eta-(\vp(\theta_{k+1} )-\vp(\theta_k ))\right].
\end{align*}
(\ref{Lang}), (\ref{generator}) and the Taylor expansion yield that
\begin{align*}
&\mathcal{A}\vp(\theta_k )\eta-(\vp(\theta_{k+1} )-\vp(\theta_k ))\\
=&\frac\eta2\Ll\nabla^2\vp(\theta_k ),\sigma\sigma^{\mathrm{T}}\Rr_\mathrm{HS}-\sqrt\eta\Ll \nabla \vp(\theta_k ),\sigma\xi_{k+1}\Rr
-\frac12\Ll\nabla^2\vp(\theta_k ),(\Delta\theta_k )(\Delta\theta_k )^\mathrm{T}\Rr_\mathrm{HS}\\
&-\frac16\int_0^1\sum_{i_1,i_1,i_3=1}^d\nabla^3_{i_1,i_2,i_3}\vp(\theta_k+t\Delta\theta_k )(\Delta \theta_k )_{i_1}(\Delta \theta_k )_{i_2}(\Delta \theta_k )_{i_3} \dif t.
\end{align*}
This, together with the previous two relations and $\Delta\theta_k=\eta g(\theta_{k})+\sqrt\eta\sigma \xi_{k+1}$, implies
\Be  \label{e:Decompose}
\eta^{-\frac 12}(\Pi_\eta(h)-\pi(h))= \mathcal{H}_\eta+\mathcal{R}_\eta,
\Ee
where, as we shall see below, $\mathcal{H}_\eta$ is a martingale and $\mathcal{R}_\eta$ is a remainder, given by
\begin{align*}
\mathcal{H}_\eta = -\eta\sum_{k=0}^{m-1}\Ll\nabla \vp(\theta_k ),\sigma\xi_{k+1}\Rr,\quad \mathcal{R}_\eta=-\sum_{i=1}^6\mathcal{R}_{\eta,i},
\end{align*}
with
\begin{align*}
\mathcal{R}_{\eta,1}=&\sqrt{\eta}(\vp(\theta_0)-\vp(\theta_m )),\
\mathcal{R}_{\eta,2}=\frac{\eta^{\frac 32}}{2}\sum_{k=0}^{m-1}\Ll\nabla^2 \vp(\theta_k ),(\sigma\xi_{k+1})(\sigma\xi_{k+1})^\mathrm{T}-\sigma\sigma^\mathrm{T}\Rr_\mathrm{HS},\\
\mathcal{R}_{\eta,3}=&\frac{\eta^2}{2}\sum_{k=0}^{m-1}\left[\Ll\nabla^2 \vp(\theta_k ),g(\theta_k )(\sigma\xi_{k+1})^\mathrm{T}\Rr_\mathrm{HS}+\Ll\nabla^2 \vp(\theta_k ),\sigma\xi_{k+1}(g(\theta_k ))^\mathrm{T}\Rr_\mathrm{HS}\right],\\
\mathcal{R}_{\eta,4}=&\frac{\eta^2}{6}\sum_{k=0}^{m-1}\int_0^1\sum_{i_1,i_2,i_3=1}^d\nabla^3_{i_1,i_2,i_3}\vp(\theta_k+t\Delta\theta_k) (\sigma\xi_{k+1})_{i_1}(\sigma\xi_{k+1})_{i_2}(\sigma\xi_{k+1})_{i_3} \dif t,\\
\mathcal{R}_{\eta,5}=&\frac{\eta^{\frac 52}}{2}\sum_{k=0}^{m-1}\Ll\nabla^2
\vp(\theta_k ),g(\theta_k )g(\theta_k )^\mathrm{T}\Rr_\mathrm{HS}\\
&+\frac{\eta^{\frac 72}}{6}\sum_{k=0}^{m-1}\int_0^1\sum_{i_1,i_2,i_3=1}^d \nabla^3_{i_1,i_2,i_3}\vp(\theta_k+t\Delta\theta_k )(g(\theta_k ))_{i_1}(g(\theta_k ))_{i_2}(g(\theta_k ))_{i_3}\dif t,\\
\mathcal{R}_{\eta,6}=&\frac{\eta^{\frac52}}{2}\sum_{k=0}^{m-1}\int_0^1\sum_{i_1,i_2,i_3=1}^d  \left[\nabla^3_{i_1,i_2,i_3}\vp(\theta_k+t\Delta\theta_k )(g(\theta_k ))_{i_1}(\sigma\xi_{k+1})_{i_2}(\sigma\xi_{k+1})_{i_3}\right.\\
&\left.+\sqrt{\eta}\nabla^3_{i_1,i_2,i_3}\vp(\theta_k+t\Delta\theta_k )(g(\theta_k ))_{i_1}(g(\theta_k ))_{i_2}(\sigma\xi_{k+1})_{i_3}
\right]\dif t.
\end{align*}

A crucial lemma for estimating the remainder $\mcl R_\eta$ is Lemma \ref{lem-g}, which has a long proof as below. To better understand the proof's strategy, we give a continuous version as the following.
For the solution $X_t$ of SDE (\ref{SDE}) and a constant $\gamma>0$ which will be chosen later, It\^{o}'s formula implies
\begin{eqnarray*}
|X_t|^2-|x|^2&=&\int_0^t2\Ll X_s,g(X_s)\Rr\dif s+\int_0^t2\Ll X_s,\sigma \dif B_s\Rr+t\|\sigma\|^2\\
&\le& -\int_0^t K_1|X_s|^2\dif s+\int_0^t2\Ll X_s,\sigma \dif B_s\Rr+(C+\|\sigma\|^2)t,
\end{eqnarray*}
where the second line follows (\ref{assu4}). Then we have
\begin{eqnarray*}
\E\exp\left\{\gamma |X_t|^2+\int_0^t \gamma K_1|X_s|^2\dif s\right\}
\le e^{\gamma|x|^2}e^{\gamma(C+\|\sigma\|^2)t}\E\exp\left\{\int_0^t2\gamma\Ll X_s,\sigma \dif B_s\Rr\right\}.
\end{eqnarray*}
H\"{o}lder's inequality and the exponential martingale property yield
\begin{eqnarray*}
&&\E\exp\left\{\int_0^t2\gamma\Ll X_s,\sigma \dif B_s\Rr\right\}\\
&\le& \left(\E\exp\left\{\int_0^t4\gamma\Ll X_s,\sigma \dif B_s\Rr-\int_0^t8\gamma^2|X_s^{\mathrm{T}}\sigma|^2\dif s\right\}\right)^{\frac12}
\left(\E\exp\left\{\int_0^t8\gamma^2|X_s^{\mathrm{T}}\sigma|^2\dif s\right\}\right)^{\frac12}\\
&=&\left(\E\exp\left\{\int_0^t8\gamma^2|X_s^{\mathrm{T}}\sigma|^2\dif s\right\}\right)^{\frac12}
\le\left(\E\exp\left\{\int_0^t\gamma K_1|X_s|^2\dif s\right\}\right)^{\frac12},
\end{eqnarray*}
where we choose $\gamma$ small enough such that $8\gamma\|\sigma\|^2\le K_1$ in the last inequality. That is
\begin{eqnarray*}
\E\exp\left\{\gamma |X_t|^2+\int_0^t \gamma K_1|X_s|^2\dif s\right\}
&\le & e^{\gamma|x|^2}e^{\gamma(C+\|\sigma\|^2)t}\left(\E\exp\left\{\int_0^t \gamma K_1|X_s|^2\dif s\right\}\right)^{\frac12} \\
&\le & e^{\gamma|x|^2}e^{\gamma(C+\|\sigma\|^2)t}\left(\E\exp\left\{\gamma |X_t|^2+\int_0^t \gamma K_1|X_s|^2\dif s\right\}\right)^{\frac12}.
\end{eqnarray*}
Hence by (\ref{assu3}), we have
\begin{eqnarray*}
\E\exp\left\{\gamma |X_t|^2+\int_0^t \frac{\gamma K_1}{2L^2}|g(X_s)|^2\dif s\right\}&\le& \E\exp\left\{\gamma |X_t|^2+\int_0^t \gamma K_1|X_s|^2\dif s\right\} e^{\frac{\gamma K_1}{L^2}|g(0)|^2t} \\
&\le& Ce^{ct}.
\end{eqnarray*}
Replacing $\frac{\gamma K_1}{2L^2}$ by $\tilde\gamma$, we can get
\begin{eqnarray*}
\E\exp\left\{\gamma |X_t|^2\right\}\le Ce^{ct}, \ \ \ \E\exp\left\{\int_0^t \tilde\gamma|g(X_s)|^2\dif s\right\}\le Ce^{ct}.
\end{eqnarray*}

\subsection{Auxiliary lemmas for $\mcl R_\eta$}
We will give in this subsection several lemmas of $\mcl R_\eta$ which play a crucial role in proving main results. Their proofs will be given in Appendix B. In order to estimate the tail probability of $\mathcal{R}_\eta$, we need the following four lemmas, the first three lemmas paving a way for proving the last.
\begin{lemma}\label{C2}
Let $\Psi_1:\R^{d}\to\R^d$ and $\Psi_2:\R^{2d}\to\R$ both be measurable functions. We have
\begin{eqnarray*}
\E_k\left[\exp\left\{\Ll \Psi_1(\theta_k),\sigma\xi_{k+1}   \Rr+\Psi_2(\theta_k,\xi_{k+1})\right\}\right]
\le \left(\E_k \left[\exp\left\{2|\Psi_1(\theta_k)|^2\|\sigma\|^2+2\Psi_2(\theta_k,\xi_{k+1})\right\}\right]\right)^{\frac 12}
\end{eqnarray*}
for $k=0,...,m-1$. Moreover, we have
\begin{eqnarray*}
\E \exp\{\sum_{k=0}^{m-1}\left(\Ll \Psi_1(\theta_k),\sigma\xi_{k+1}   \Rr+\Psi_2(\theta_k,\xi_{k+1})\right)\}
\le \left(\E \exp\{\sum_{k=0}^{m-1}2\left(|\Psi_1(\theta_k)|^2\|\sigma\|^2+\Psi_2(\theta_k,\xi_{k+1})\right)\}\right)^{\frac 12},
\end{eqnarray*}
and
\begin{eqnarray*}
\E_0 \exp\{\sum_{k=0}^{m-1}\left(\Ll \Psi_1(\theta_k),\sigma\xi_{k+1}   \Rr+\Psi_2(\theta_k,\xi_{k+1})\right)\}
\le \left(\E_0 \exp\{\sum_{k=0}^{m-1}2\left(|\Psi_1(\theta_k)|^2\|\sigma\|^2+\Psi_2(\theta_k,\xi_{k+1})\right)\}\right)^{\frac 12}.
\end{eqnarray*}
\end{lemma}

\begin{lemma}\label{lem-g}
Under Assumption \ref{assu}, there exist $\eta_0>0$ and $\gamma_0>0$, both depending on $L,K_1,K_2$,$|g(0)|^2$ and $\sigma$,  such that as $\eta<\eta_0$ and $\gamma<\gamma_0$,
\begin{align}\label{R_1,22-0}
\E_0\exp\left\{\gamma \eta\sum_{k=0}^{m-1}|g(\theta_k)|^2\right\}\le C e^{c (\eta^{-1}+|\theta_0|^2)},
\end{align}
where $C$ and $c$ depend on $L,K_1,K_2$,$|g(0)|^2$, $\sigma$ and $\gamma$.
Moreover, if $\theta_0 \sim \pi_\eta$,
\begin{align}\label{R_1,22}
\E\exp\left\{\gamma \eta\sum_{k=0}^{m-1}|g(\theta_k)|^2\right\}\le C e^{c \eta^{-1}},
\end{align}
where $C$ and $c$ depend on $L,K_1,K_2$,$|g(0)|^2$, $\sigma$ and $\gamma$. This particular implies that for all $x>0$,
\begin{align}\label{e:R-1-22-0}
\PP_0\left(\eta\sum_{k=0}^{m-1}|g(\theta_k)|^2>x\right)\le C e^{c_1 (\eta^{-1}+|\theta_0|^2)} e^{-c_2 x},
\end{align}
where $C$, $c_1$, $c_2$ depends on $L,K_1,K_2$,$|g(0)|^2$, $\sigma$ and $\gamma_0$. Moreover, if $\theta_0 \sim \pi_\eta$,
\begin{align}\label{e:R-1-22}
\PP\left(\eta\sum_{k=0}^{m-1}|g(\theta_k)|^2>x\right)\le C e^{c_1 \eta^{-1}} e^{-c_2 x},
\end{align}
where $C$, $c_1$, $c_2$ depends on $L,K_1,K_2$,$|g(0)|^2$, $\sigma$ and $\gamma_0$.

\end{lemma}


\begin{lemma}\label{C5}
Let $\Psi:\R^{2d}\to\R^d$ be measurable function satisfying the conditions
$$
\E_k[\Psi(\theta_k,\xi_{k+1})]=0\quad\text{and}\quad \PP_k(|\Psi(\theta_k,\xi_{k+1})|\le K(1+|\xi_{k+1}|^2))=1
$$
for $k=0,...,m-1$, where $K \in (0,\infty)$ is an arbitrary constant. Then we have
$$
\E\left[\exp\left\{\frac{1}{\sqrt m}\sum_{n=0}^{m-1}\Psi(\theta_k,\xi_{k+1})\right\}\right]\le C,
$$
where $C$ depends on $K$.
\end{lemma}

\begin{lemma}\label{lem-R}
Suppose that Assumption \ref{assu} holds. Let $h \in C^2_b(\R^d,\R)$ and $\vp:\R^d\to\R$ be the solution of (\ref{stein}). We have
\begin{eqnarray*}
\PP_0(|\mathcal{R}_\eta|>x)\le& Ce^{c_1\|\theta_0\|^2} \left(e^{-c_2\eta^{-\frac12}x^{\frac12} }1_{\{x<\eta^{-1}\}}+ e^{-c_2\eta^{-\frac 35}x^{\frac 25}}1_{\{x\ge \eta^{-1}\}}+e^{-c_2\eta^{-2\bar\gamma}x^{\frac23}}\right),
\end{eqnarray*}
where $0<\bar\gamma<\frac{1}{4}$  and $x\ge c\max\{\eta^{\frac32-6\bar\gamma},\eta^{\frac32\bar\gamma},\eta^{\frac 12}\}$. Here $C$, $c$, $c_1$, $c_2$ depends on $L,K_1,K_2$,$|g(0)|^2$, $\sigma$. Moreover, for $\theta_0\sim\pi_\eta$, we have
\begin{eqnarray*}
\PP(|\mathcal{R}_\eta|>x)\le& C \left(e^{-c_1\eta^{-\frac12}x^{\frac12} }1_{\{x<\eta^{-1}\}}+ e^{-c_1\eta^{-\frac 35}x^{\frac 25}}1_{\{x\ge \eta^{-1}\}}+e^{-c_1\eta^{-2\bar\gamma}x^{\frac23}}\right),
\end{eqnarray*}
where $0<\bar\gamma<\frac{1}{4}$  and $x\ge c\max\{\eta^{\frac32-6\bar\gamma},\eta^{\frac32\bar\gamma},\eta^{\frac 12}\}$. Here $C$, $c$, $c_1$ depends on $L,K_1,K_2$,$|g(0)|^2$, $\sigma$.
\end{lemma}

\section{Proof of Theorem \ref{thm-CLT}}
We first introduce following lemma which paves a way to proving the convergence of martingale $\mathcal{H}_\eta$. Its proof borrows the idea of the Stein's method in \citet[Theorem 2.5]{FSX19}.
\begin{lemma}\label{lu1}
Let $\pi$ and $\pi_\eta$ be the same as those in Lemma \ref{l:Erg2}, $\vp$ be the solution of Stein's equation \eqref{stein}. We have
\begin{align*}
|\pi_\eta(|\sigma^{\rm T}\nabla\vp|^2)-\pi(|\sigma^{\rm T}\nabla\vp|^2)|\le C\eta^{\frac12}.
\end{align*}
\end{lemma}
Here, $C$ depends on $g$ and $\sigma$.
\begin{proof}
We shall use the stationary Markov chain trick in \citet[Theorem 2.5]{FSX19}.
Let $\{\theta_k \}_{k\ge0}$ be the Markov chain with initial value $\theta_0 \sim\pi_\eta$. (\ref{Lang}) implies that
\begin{align}\label{e:t2.4-3}
&\E_0[\Delta \theta_0 ]=\eta g(\theta_0 ),\quad\E_0[(\Delta \theta_0 )(\Delta \theta_0 )^{\mathrm{T}}]=\eta^2 g(\theta_0 )g^{\mathrm{T}}(\theta_0 )+\eta\sigma\sigma^{\mathrm{T}}.
\end{align}
{For Stein's equation
\begin{align}\label{stein1}
|\sigma^{\rm T}\nabla\vp|^2-\pi(|\sigma^{\rm T}\nabla\vp|^2)=\mathcal A{\bar\vp},
\end{align}
Lemma \ref{solpoisson} implies the test function $|\sigma^{\rm T}\nabla\vp|^2\in C^2_b(\R^d,\R)$. Thus $\bar\vp$ exists and satisfies $\|\nabla^k \bar\vp\|\le C$ with $k=0,1,2,3,4$ by Lemma \ref{solpoisson} again.
}
The Taylor expansion and the stationarity of $(\theta_k)_{k \ge 0}$ yield
\begin{align}\label{xu1}
0=&\E[\bar\vp(\theta_1)-\bar\vp(\theta_0)]\\
=&\E[\Ll\nabla \bar\vp(\theta_0),\Delta\theta_0 \Rr]+\frac{1}{2}\E[\Ll\nabla^2 \bar\vp(\theta_0),\Delta\theta_0(\Delta\theta_0)^{\mathrm{T}}\Rr_{\mathrm{HS}}]\nonumber\\
&+\frac16\int_0^1\E\left[\sum_{i_1,i_1,i_3=1}^d\nabla^3_{i_1,i_2,i_3}\bar\vp(\theta_0+t\Delta\theta_0 )(\Delta\theta_0 )_{i_1}(\Delta\theta_0 )_{i_2}(\Delta\theta_0 )_{i_3}\right]\dif t.\nonumber\nonumber
\end{align}
For the first and the second terms, by (\ref{e:t2.4-3}), we have
\begin{align*}
\E[\Ll\nabla \bar\vp(\theta_0),\Delta\theta_0 \Rr]=\E[\Ll\nabla \bar\vp(\theta_0),\E_0[\Delta\theta_0 ]\Rr]=
\E[\Ll\nabla \bar\vp(\theta_0),\eta g(\theta_0 )\Rr],
\end{align*}
\begin{align*}
\E[\Ll\nabla^2 \bar\vp(\theta_0),\Delta\theta_0(\Delta\theta_0)^{\mathrm{T}}\Rr_{\mathrm{HS}}]&=
\E[\Ll\nabla^2 \bar\vp(\theta_0),\E_0[\Delta\theta_0(\Delta\theta_0)^{\mathrm{T}}]\Rr_{\mathrm{HS}}]\\
&=
\E[\Ll\nabla^2 \bar\vp(\theta_0),\eta^2 g(\theta_0 )g^{\mathrm{T}}(\theta_0 )+\eta\sigma\sigma^{\mathrm{T}}\Rr_{\mathrm{HS}}].
\end{align*}
Combining equalities above with (\ref{generator}) and (\ref{xu1}), we have
\begin{eqnarray*}
\E[\mathcal{A}(\bar\vp(\theta_0))]&=&-\frac{1}{2}\E[\Ll\nabla^2\bar\vp(\theta_0),\eta g(\theta_0 )g^{\mathrm{T}}(\theta_0 )\Rr_{\rm {HS}}]\\
&&-\frac{1}{6\eta}\int_0^1\E\left[\sum_{i_1,i_1,i_3=1}^d\nabla^3_{i_1,i_2,i_3}\bar\vp(\theta_0+t\Delta\theta_0)(\Delta\theta_0 )_{i_1}
(\Delta\theta_0 )_{i_2}(\Delta\theta_0 )_{i_3}\right]\dif t.
\end{eqnarray*}
For the first term, the boundedness of $\|\nabla^2\bar\vp\|$ and (\ref{e:BenDan-1}) imply
\begin{eqnarray*}
\left|\frac{1}{2}\E[\Ll\nabla^2\bar\vp(\theta_0),\eta g(\theta_0 )g^{\mathrm{T}}(\theta_0 )\Rr_{\rm {HS}}]\right|\le C\pi_\eta(|g|^4)^{1/2}\eta\le C\eta.
\end{eqnarray*}
For the second term,  by (\ref{e:BenDan-1}), we can get
\begin{align*}
\E[|g(\theta_0)|^3]\le\big(\E[|g(\theta_0)|^4]\big)^{3/4}=\pi_\eta(|g|^4)^{3/4}<\infty.
\end{align*}
Cauchy's inequality and the boundedness of $\|\nabla^3\bar\vp\|$ imply
\begin{align*}
&\left|\frac{1}{6\eta}\int_0^1\E\left[\sum_{i_1,i_1,i_3=1}^d\nabla^3_{i_1,i_2,i_3}\bar\vp(\theta_0+t\Delta\theta_0)(\Delta\theta_0 )_{i_1}
(\Delta\theta_0 )_{i_2}(\Delta\theta_0 )_{i_3}\right]\dif t\right|\\
\le& \frac{C}{\eta}\int_0^1\E[\|\nabla^3\bar\vp(\theta_0+t\Delta\theta_0)\||\Delta(\theta_0)|^3] \dif t
\le C(\eta^2\E[|g(\theta_0 )|^3]+\eta^{\frac12}\E[|\xi_1|^3])
\le C\eta^{\frac12}.
\end{align*}
Here, the constant $C$ depends on $\sigma$ and $g$. Hence we have
$$|\E[\mathcal{A}(\bar\vp(\theta_0))]| \le C\eta^{\frac12}.$$

From Stein's equation (\ref{stein1}), we deduce
\begin{align*}
|\pi_\eta(|\sigma^{\rm T}\nabla\vp|^2)-\pi(|\sigma^{\rm T}\nabla\vp|^2)|=|\E[|\sigma^{\rm T}\nabla\vp(\theta_0)|^2-\pi(|\sigma^{\rm T}\nabla\vp|^2)]|=|\E[\mathcal{A}(\bar\vp(\theta_0))]|
\le C\eta^{\frac12}.
\end{align*}
\end{proof}

{\begin{lemma}\label{l:clt}
Under the condition of Theorem \ref{thm-CLT}, we have
$$
\mathcal{H}_\eta\Rightarrow N(0, \pi(|\sigma^{\mathrm{T}}\nabla\vp|^2)).
$$
\end{lemma}
\begin{proof}
Recall $\mathcal{H}_\eta=-\eta\sum_{i=0}^{m-1}\Ll\nabla\vp(\theta_i),\sigma\xi_{i+1}\Rr$. We denote
$$
Z_i=\Ll\nabla\vp(\theta_i),\sigma\xi_{i+1}\Rr,\quad i=0,...,m-1.$$
\citet[Theorem 2.3]{McL1974} will imply the result if we can verify the conditions
\begin{eqnarray}
&&\E\max_{0\le i\le m-1}\{\eta|Z_i|\}\to0,\label{e:cond1}\\
&&\eta^2\sum_{i=0}^{m-1}Z_i^2\to \pi(|\sigma^{\mathrm{T}}\nabla \vp|^2)\text{~in probability.}\label{e:cond2}
\end{eqnarray}
Denoting $\hat{Z_i}=Z_i1_{\{|Z_i|^2\le \eta^{-1}\}}$ and
$\check{Z_i}=Z_i1_{\{|Z_i|^2> \eta^{-1}\}}$, we have
\begin{eqnarray*}
\eta^2(\max_{0\le i\le m-1}\{|Z_i|\})^2=\eta^2\max_{0\le i\le m-1}\{|Z_i|^2\}\le \eta^2\max_{0\le i\le m-1}\{|\hat{Z_i}|^2\}+\eta^2\max_{0\le i\le m-1}\{|\check{Z_i}|^2\}.
\end{eqnarray*}
It is easily to see that the first term converges to $0$ in probability. For the second term, we have
\begin{eqnarray*}
\eta^2\E\max_{0\le i\le m-1}\{|\check{Z_i}|^2\}&\le &\eta^2\sum_{0\le i\le m-1} \E |\check{Z_i}|^2.
\end{eqnarray*}
Since $\E[Z_i^2]$ is finite, $ \E|\check Z_i|^2 $ converges to 0 as $\eta\to0$ for each $i$, this implies that
$\eta^2\E[\max_{0\le i\le m-1}\{|\check{Z_i}|^2\}]$ converges to $0$. H\"{o}lder's inequality yields (\ref{e:cond1}).

For (\ref{e:cond2}), we can finish the proof if we verify
\begin{eqnarray}\label{e:3.7}
&&\E[\eta^2\sum_{i=0}^{m-1}\left(Z_i^2-\pi(|\sigma^{\mathrm{T}}\nabla\vp|^2)\right)]^2\\
&\le&2\E[\eta^2\sum_{i=0}^{m-1}\left(Z_i^2-\pi_\eta(|\sigma^{\mathrm{T}}\nabla\vp|^2)\right)]^2
+2\Big(\eta^2\sum_{i=0}^{m-1}\big(\pi_\eta(|\sigma^{\mathrm{T}}\nabla\vp|^2-\pi(|\sigma^{\mathrm{T}}\nabla\vp|^2)\big)\Big)^2\to0.\nonumber
\end{eqnarray}
By Lemma \ref{lu1}, the second term converges to $0$.  For the first term, a straight calculation gives that
\begin{eqnarray*}
\E[\eta^2\sum_{i=0}^{m-1}\left(Z_i^2-\pi_\eta(|\sigma^{\mathrm{T}}\vp|^2)\right)]^2
&=&\eta^4\sum_{i=0}^{m-1}\E[Z_i^2-\pi_\eta(|\sigma^{\mathrm{T}}\vp|^2)]^2\\
&&+2\eta^4\sum_{i,j=0,i<j}^{m-1}\E\left[(Z_i^2-\pi_\eta(|\sigma^{\mathrm{T}}\vp|^2))(Z_j^2-\pi_\eta(|\sigma^{\mathrm{T}}\vp|^2))\right].
\end{eqnarray*}
For the first term, the boundedness of $\|\nabla\vp\|$ implies
\begin{eqnarray}\label{e:3.4}
\E[Z_i^2-\pi_\eta(|\sigma^{\mathrm{T}}\vp|^2)]^2
&\le& 2\E[Z_i^4]+2\pi_\eta(|\sigma^{\mathrm{T}}\vp|^2)^2\nonumber\\
&\le& 2\E[|\nabla\vp(\theta_i)|^4\|\sigma\|^4|\xi_{i+1}|^4] +2\pi_\eta(|\sigma^{\mathrm{T}}\vp|^2)^2\nonumber\\
&\le& C+2\pi_\eta(|\sigma^{\mathrm{T}}\vp|^2)^2.
\end{eqnarray}
Then we have
\begin{eqnarray*}
\eta^4\sum_{i=0}^{m-1}\E[Z_i^2-\pi_\eta(|\sigma^{\mathrm{T}}\vp|^2)]^2&\le& C\eta^4 m \to0.
\end{eqnarray*}

For the second term, we can calculate that
\begin{eqnarray*}
&&\sum_{i,j=0,i<j}^{m-1}\E\left[(Z_i^2-\pi_\eta(|\sigma^{\mathrm{T}}\vp|^2))(Z_j^2-\pi_\eta(|\sigma^{\mathrm{T}}\vp|^2))\right]\\
&=&\sum_{i,j=0,i<j}^{m-1}\E\left[(Z_i^2-\pi_\eta(|\sigma^{\mathrm{T}}\vp|^2))\E_{i+1}[Z_j^2-\pi_\eta(|\sigma^{\mathrm{T}}\vp|^2)]\right]\\
&=&\sum_{i,j=0,i<j}^{m-1}\E\left[(Z_i^2-\pi_\eta(|\sigma^{\mathrm{T}}\vp|^2))\E_{i+1}[|\sigma^{\rm T}\nabla\vp(\theta_j)|^2-\pi_\eta(|\sigma^{\mathrm{T}}\vp|^2)]\right]\\
&\le&\sum_{i,j=0,i<j}^{[t]}\E\left[|Z_i^2-\pi_\eta(|\sigma^{\mathrm{T}}\vp|^2)|(1+|\theta_{i+1}|^2) e^{-c(j-i-1)} \right],
\end{eqnarray*}
where the last inequality follows from (\ref{ergodic2}). By H\"{o}lder's inequality, we have
\begin{eqnarray*}
&&\sum_{i,j=0, i<j}^{m-1}\E\left[(Z_i^2-\pi_\eta(|\sigma^{\mathrm{T}}\vp|^2))(Z_j^2-\pi_\eta(|\sigma^{\mathrm{T}}\vp|^2))\right]\\
&\le&\sum_{i,j=0, i<j}^{m-1}{\left[\left(\E\left[(Z_i^2-\pi_\eta(|\sigma^{\mathrm{T}}\vp|^2))\right]^2\right)^{\frac12}\left(\E\left[(1+|\theta_{i+1}|^2) e^{-c(j-i-1)}\right]^2\right)^{\frac12}\right]}\\
&\le& C\sum_{i,j=0, i<j}^{m-1}e^{-c(j-i)}\left(1+\E|\theta_{i+1}|^4\right)^{\frac12},
\end{eqnarray*}
where the boundedness of $\E\left[(Z_i^2-\pi_\eta(|\sigma^{\mathrm{T}}\vp|^2))\right]^2$ follows from (\ref{e:3.4}).
Now we estimate $\E|\theta_{i+1}|^4$. A similar calculation with (\ref{e:BenDan-0}) yields
\begin{eqnarray*}
\E[|\theta_{i+1}|^4]=\E[\E_i|\theta_{i+1}|^4]\le (1-K_1\eta+c\eta^2)\E[|\theta_i|^4]+C\eta.
\end{eqnarray*}
By iteration with initial data $\theta_0=x$, we obtain
\begin{eqnarray*}
\E[|\theta_{i+1}|^4]\le C\eta\sum_{k=0}^{i}(1-K_1\eta+c\eta^2)^k+|x|^4(1-K_1\eta+c\eta^2)^{i+1}.
\end{eqnarray*}
Choosing $\eta$ small enough such that $1-K_1\eta+c\eta^2<1$ gives
\begin{eqnarray*}
\E[|\theta_{k}|^4]\le |x|^4+\frac{C}{K_1+c\eta},\quad k=0,1,...
\end{eqnarray*}

Combining the relationships above, we have
\begin{eqnarray*}
&&2\eta^4\sum_{i,j=0,i<j}^{m-1}\E\left[(Z_i^2-\pi_\eta(|\sigma^{\mathrm{T}}\vp|^2))(Z_j^2-\pi_\eta(|\sigma^{\mathrm{T}}\vp|^2))\right]
\le C\eta^4\sum_{i,j=0,\atop i<j}^{m-1}e^{-c(j-i)}\\
&=&C\eta^4\sum_{i,j=0,\atop 0<j-i\le \ln m}^{m-1}e^{-c(j-i)}+C\eta^4\sum_{i,j=0,\atop \ln m<j-i }^{m-1}e^{-c(j-i)}\\
&\le& C\eta^4 m\ln m+C\eta^4 e^{-c\ln m}(m-\ln m)^2\to0.
\end{eqnarray*}
Hence we prove the first term of \eqref{e:3.7} converges to $0$ and finish the proof.

\end{proof}
}

\begin{proof}[Proof of Theorem \ref{thm-CLT}]
We have shown in \eqref{e:Decompose} that
\begin{eqnarray*}
\eta^{-\frac 12}(\Pi_\eta(h)-\pi(h))= \mathcal{H}_\eta+\mathcal{R}_\eta.
\end{eqnarray*}
Here $\mathcal{H}_\eta$ weakly converges to $N(0,\pi(|\sigma^{\rm T}\nabla\vp|^2))$ by Lemma \ref{l:clt}. Lemma \ref{lem-R} implies $\mathcal{R}_\eta$ converges to $0$ in probability with fixed initial value $\theta_0$. Thus $\eta^{-\frac 12}(\Pi_\eta(h)-\pi(h))\Rightarrow N(0,\pi(|\sigma^{\rm T}\nabla\vp|^2))$.

\end{proof}

\section{Proof of Theorem \ref{thm-self-normal}}
\subsection{Self-normalized Cram\'er-type moderate deviation of $\mcl H_\eta$}

In order to prove the Cram\'er-type moderate deviation result for $\mcl H_\eta$, we introduce following concentration inequality for stationary process.
\begin{lemma}\label{EY-Y}
Suppose that the conditions of Theorem \ref{thm-self-normal} hold. Then, for any $y>0$
\begin{eqnarray*}
\PP\left(\left|\sum_{i=0}^{k-1}|\sigma^{\rm{T}}\vp(\theta_i)|^2-k\pi_\eta(|\sigma^{\rm{T}}\vp|^2)\right|>y\right)
\le 2e^{-Cy^2k^{-1}},\quad k\in\mathbb{N}.
\end{eqnarray*}
\end{lemma}
Here, $C$ depends on $g$ and $\sigma$.
\begin{proof}
Since $\theta_0\sim\pi_\eta$, $(\theta_k )_{k\ge0}$ is stationary. Following \citet[(6)]{dedecker2015subgaussian} with
$||\sigma^{\rm{T}}\vp(\theta_k)|^2-\pi_\eta(|\sigma^{\rm{T}}\vp|^2)|\le C$, we can get the result immediately.
\end{proof}

\begin{lemma}\label{Fan_Coro2.1}
Under the conditions of Theorem \ref{thm-self-normal}, one has
\begin{align*}
\PP\left( \frac{\mathcal{H}_\eta}{\sqrt{ \Y_\eta }} \geq x\right) \Big/(1-\Phi(x))  = 1+O(x\eta^{\frac13}+\eta^{\frac13}),
\end{align*}
uniformly for $\eta^{\frac13}\le x=o(\eta^{-\frac 13})$ as $\eta$ tends to zero. Here, $O$ and $o$ depend on $g$, $\sigma$.
\end{lemma}
\begin{proof}
We first prove the upper bound of $\PP\left( \frac{\mathcal{H}_\eta}{\sqrt{\Y_\eta }} \geq x\right)  \Big/(1-\Phi(x))$. Notice that $\E\Y_\eta=\pi_\eta(|\sigma^{\rm T}\nabla\vp|^2)$ by the fact $\theta_0\sim\pi_\eta$, without loss of generality, we may assume $\E \Y_\eta=\pi_\eta(|\sigma^{\rm T}\nabla\vp|^2)=1$. For $y$ such that $0<y\eta^2<1$ which will be chosen later, Lemma \ref{EY-Y} implies
\begin{eqnarray}\label{e:sum4}
\PP\left(\mathcal{H}_\eta/\sqrt{\Y_\eta } \geq x\right)
&=&\PP\left(\mathcal{H}_\eta/\sqrt{ \Y_\eta } \geq x, \eta^{-2}\left|1 - \Y_\eta\right|>y \right)\nonumber+\PP\left(\mathcal{H}_\eta/\sqrt{ \Y_\eta } \geq x, \eta^{-2}\left|1 - \Y_\eta\right|\le y \right)\nonumber\\
&\le&\PP\left(\eta^{-2}\left|1 -\Y_\eta\right|>y \right)+\PP\left(\mathcal{H}_\eta/\sqrt{1 -y\eta^2} \geq x, \eta^{-2}\left|1 - \Y_\eta\right|\le y \right)\nonumber\\
&\le&2e^{-Cy^2\eta^2}+\PP\left(\mathcal{H}_\eta/\sqrt{1 -y\eta^2} \geq x\right).
\end{eqnarray}
Define $$\tilde{\mathcal{H}}_t=\eta\int_0^t\nabla\vp(\theta_{[s]})^{\rm T}\sigma\dif B_s$$
for any $t\in\R^+$ which is a continuous martingale. Denote its sharp bracket process by $\Ll\tilde{\mathcal{H}}\Rr(s,t)=\eta^2\int_s^t|\sigma^{\rm T}\nabla\vp(\theta_{[r]})|^2\dif r$ and $\Ll\tilde{\mathcal{H}}\Rr(t)=\Ll\tilde{\mathcal{H}}\Rr(0,t)$ for simplicity. It is easy to see
$$\tilde{\mathcal{H}}_m\overset{d}{=}\mathcal{H}_\eta, $$
 $$\Ll\tilde{\mathcal{H}}\Rr(t)=\eta^2\int_0^t |\sigma^{\rm T}\nabla\vp(\theta_{[s]})|^2\dif s=\sum_{i=0}^{[t]-1}\eta^2|\sigma^{\rm T}\nabla\vp(\theta_{i})|^2+\eta^2\int_{[t]}^t|\sigma^{\rm T}\nabla\vp(\theta_{[s]})|^2\dif s. $$
Denoting the stopping time $T_1=\inf\{s:\Ll\tilde{\mathcal{H}}\Rr(s)>1\}$, Dambis-Dubins-Schwarz  Theorem (cf. \citet[Theorem 5.1.6]{revuz2013continuous}) yields that $\tilde{\mathcal{H}}_{T_1}$ is a $\mathcal{F}_{T_1}-$Brownian motion and $\tilde{\mathcal{H}}_{T_1}\sim N(0,1)$. Then we have
\begin{eqnarray}\label{e:sum5}
\PP\left(\mathcal{H}_\eta/\sqrt{1 -y\eta^2} \geq x\right)
&=&\PP\left(\frac{\tilde{\mathcal{H}}_m-\tilde{\mathcal{H}}_{T_1}+\tilde{\mathcal{H}}_{T_1}}{\sqrt{1-y\eta^2}} \geq x \right)\\
&\le&\PP\left(\frac{\tilde{\mathcal{H}}_m-\tilde{\mathcal{H}}_{T_1}}{\sqrt{1-y\eta^2}} \geq c_0 \right)
+\PP\left(\frac{\tilde{\mathcal{H}}_{T_1}}{\sqrt{1-y\eta^2}} \geq x-c_0 \right)\nonumber
\end{eqnarray}
with small $c_0$ satisfying $0<c_0\le x$ which will be chosen later. For the second term on the right hand side, since $\tilde{\mathcal{H}}_{T_1}\sim N(0,1)$,
\begin{eqnarray}\label{e:sum6}
\PP\left(\frac{\tilde{\mathcal{H}}_{T_1}}{\sqrt{1-y\eta^2}} \geq x-c_0 \right)=1-\Phi(\sqrt{1-y\eta^2}(x-c_0)).
\end{eqnarray}

For the first term and $\alpha\in(0,1)$, we have
\begin{eqnarray}\label{e:6}
\PP\left(\frac{\tilde{\mathcal{H}}_m-\tilde{\mathcal{H}}_{T_1}}{\sqrt{1-y\eta^2}} \geq c_0 \right)
&=&\PP\left(\frac{\tilde{\mathcal{H}}_m-\tilde{\mathcal{H}}_{T_1}}{\sqrt{1-y\eta^2}} \geq c_0, T_1\notin [m-m^\alpha,m+m^\alpha] \right)\\
&&+\PP\left(\frac{\tilde{\mathcal{H}}_m-\tilde{\mathcal{H}}_{T_1}}{\sqrt{1-y\eta^2}} \geq c_0, T_1\in [m-m^\alpha,m+m^\alpha] \right).\nonumber
\end{eqnarray}
Without loss of generality, we may assume that $m^\alpha$ is an integer. The definition of $T_1$ implies $\{T_1< m-m^\alpha\}=\{\Ll\tilde{\mathcal{H}}\Rr(m-m^\alpha)>1\}$ and $\{T_1> m+m^\alpha\}=\{\Ll\tilde{\mathcal{H}}\Rr(m+m^\alpha)<1\}$. Then we can obtain
\begin{eqnarray*}
\PP\big(\frac{\tilde{\mathcal{H}}_m-\tilde{\mathcal{H}}_{T_1}}{\sqrt{1-y\eta^2}} \geq c_0, T_1\notin [m-m^\alpha,m+m^\alpha] \big)
&\le&\PP\left(T_1< m-m^\alpha\right)
+\PP\left(T_1> m+m^\alpha\right)\\
&=&\PP\left(\Ll\tilde{\mathcal{H}}\Rr(m-m^\alpha)>1\right)
+\PP\left(\Ll\tilde{\mathcal{H}}\Rr(m+m^\alpha)<1\right).
\end{eqnarray*}
Following Lemma \ref{EY-Y},  one has
\begin{eqnarray*}
\PP\left(\Ll\tilde{\mathcal{H}}\Rr(m-m^\alpha)>1\right)
&=&\PP\left(    \sum_{i=0}^{m-m^\alpha-1}|\sigma^{\rm T}\nabla\vp(\theta_i)|^2-(m-m^\alpha)>m-(m-m^\alpha)\right)\\
&\le& e^{-Cm^{2\alpha}(m-m^\alpha)^{-1}}\le e^{-Cm^{2\alpha-1}}.
\end{eqnarray*}
Similarly, we can get
$\PP\left(\Ll\tilde{\mathcal{H}}\Rr(m+m^\alpha)<1\right)
\le e^{-Cm^{2\alpha-1}}.$ That is
\begin{eqnarray}\label{e:sum7}
\PP\left(\frac{\tilde{\mathcal{H}}_m-\tilde{\mathcal{H}}_{T_1}}{\sqrt{1-y\eta^2}} \geq c_0, T_1\notin [m-m^\alpha,m+m^\alpha] \right)\le 2e^{-Cm^{2\alpha-1}}.
\end{eqnarray}
For the second term of \eqref{e:6}, we have
\begin{eqnarray}\label{e:sum10}
&&\PP\left(\frac{\tilde{\mathcal{H}}_m-\tilde{\mathcal{H}}_{T_1}}{\sqrt{1-y\eta^2}} \geq c_0, T_1\in [m-m^\alpha,m+m^\alpha] \right)\\
&\le&\PP\left(\sup_{s\in[m-m^\alpha,m+m^\alpha]}(\tilde{\mathcal{H}}_m-\tilde{\mathcal{H}}_{s})  \geq c_0\sqrt{1-y\eta^2} \right)\nonumber\\
&\le&\PP\left(\sup_{s\in[m,m+m^\alpha]}(\tilde{\mathcal{H}}_m-\tilde{\mathcal{H}}_{s})  \geq c_0\sqrt{1-y\eta^2} \right)
+\PP\left(\sup_{s\in[m-m^\alpha,m]}(\tilde{\mathcal{H}}_m-\tilde{\mathcal{H}}_{s})  \geq c_0\sqrt{1-y\eta^2} \right).\nonumber
\end{eqnarray}
For the first term and positive number $y'$ which will be chosen later, one has
\begin{eqnarray*}
&&\PP\left(\sup_{s\in[m,m+m^\alpha]}(\tilde{\mathcal{H}}_m-\tilde{\mathcal{H}}_s)  \geq c_0\sqrt{1-y\eta^2} \right)\\
&\le&\PP\left(\Ll \tilde{\mathcal{H}}\Rr(m,m+m^\alpha)>\eta^2y'\right)\\
&&+\PP\left(\sup_{s\in[m,m+m^\alpha]}(\tilde{\mathcal{H}}_{s}-\tilde{\mathcal{H}}_m)  \geq c_0\sqrt{1-y\eta^2} , \Ll \tilde{\mathcal{H}}\Rr(m,m+m^\alpha)\le \eta^2y'\right).
\end{eqnarray*}
The first probability can be estimate by Lemma \ref{EY-Y}, that is,
\begin{eqnarray*}
\PP\left(\Ll \tilde{\mathcal{H}}\Rr(m,m+m^\alpha)>\eta^2y'\right)\le e^{-\frac{C(y'-m^\alpha)^2}{m^\alpha}}.
\end{eqnarray*}
For the second probability, the Bernstein inequality (cf. \citet[Proposition 4.2.3(1)]{barlow1986inequalities}) implies
\begin{eqnarray*}
&&\PP\left(\sup_{s\in[m,m+m^\alpha]}(\tilde{\mathcal{H}}_{s}-\tilde{\mathcal{H}}_m)  \geq c_0\sqrt{1-y\eta^2} , \Ll \tilde{\mathcal{H}}\Rr(m,m+m^\alpha)\le \eta^2y'\right)\\
&\le& \exp\{-\frac{c_0^2(1-y\eta^2)}{2\eta^2y'}\}= e^{-\frac{c_0^2(m-y)}{2y'}}.
\end{eqnarray*}
Thus we have
\begin{eqnarray}\label{e:sum9}
\PP\left(\sup_{s\in[m,m+m^\alpha]}(\tilde{\mathcal{H}}_m-\tilde{\mathcal{H}}_{s})  \geq c_0\sqrt{1-y\eta^2} \right)\le e^{-\frac{C(y'-m^\alpha)^2}{m^\alpha}}+e^{-\frac{c_0^2(m-y)}{2y'}}.
\end{eqnarray}
For the second term of \eqref{e:sum10},
\begin{eqnarray*}
&&\PP\left(\sup_{s\in[m-m^\alpha,m]}(\tilde{\mathcal{H}}_m-\tilde{\mathcal{H}}_{s})  \geq c_0\sqrt{1-y\eta^2} \right)\\
&\le&\sum_{k=0}^{m^\alpha-1}\PP\left(\sup_{s\in[m-k-1,m-k]}(\tilde{\mathcal{H}}_m-\tilde{\mathcal{H}}_{s})  \geq c_0\sqrt{1-y\eta^2} \right)\\
&\le&\sum_{k=0}^{m^\alpha-1}\PP\left(\tilde{\mathcal{H}}_m-\tilde{\mathcal{H}}_{m-k}  \geq \frac{c_0}{2}\sqrt{1-y\eta^2} \right)
+\sum_{k=0}^{m^\alpha-1}\PP\left(\sup_{s\in[m-k-1,m-k]}(\tilde{\mathcal{H}}_{m-k}-\tilde{\mathcal{H}}_{s})  \geq \frac{c_0}{2}\sqrt{1-y\eta^2} \right).
\end{eqnarray*}
For the first probability, the stability of $\theta_k$ and  \eqref{e:sum9} yield
\begin{eqnarray*}
\sum_{k=0}^{m^\alpha-1}\PP\left(\tilde{\mathcal{H}}_m-\tilde{\mathcal{H}}_{m-k}  \geq \frac{c_0}{2}\sqrt{1-y\eta^2} \right)
&=&\sum_{k=0}^{m^\alpha-1}\PP\left(\tilde{\mathcal{H}}_{m+k}-\tilde{\mathcal{H}}_{m}  \geq \frac{c_0}{2}\sqrt{1-y\eta^2} \right)
\\
&\le& m^\alpha \PP\left(\sup_{s\in[m,m+m^\alpha]}(\tilde{\mathcal{H}}_{s}-\tilde{\mathcal{H}}_{m})  \geq \frac{c_0}{2}\sqrt{1-y\eta^2} \right)\\
&\le& m^{\alpha} \left(e^{-\frac{C(y'-m^\alpha)^2}{m^\alpha}}+e^{-\frac{c_0^2(m-y)}{8y'}}\right).
\end{eqnarray*}
For the second probability, by the boundedness of $\nabla\vp$, we have
\begin{eqnarray*}
&&\sum_{k=0}^{m^\alpha-1}\PP\left(\sup_{s\in[m-k-1,m-k]}(\tilde{\mathcal{H}}_{m-k}-\tilde{\mathcal{H}}_{s})  \geq \frac{c_0}{2}\sqrt{1-y\eta^2} \right)\\
&=&\sum_{k=0}^{m^\alpha-1}\PP\left(\sup_{s\in[m-k-1,m-k]}\int_s^{m-k}\eta(\nabla\vp(\theta_{[r]}))^{\rm T}\sigma\dif Br\geq \frac{c_0}{2}\sqrt{1-y\eta^2} \right)\\
&\le&\sum_{k=0}^{m^\alpha-1}\PP\left(\sup_{m-k-1\le s\le t\le m-k}|B_t-B_s|\geq \frac{C c_0}{\eta}\sqrt{1-y\eta^2} \right)\\
&\le&m^\alpha\PP\left(\sup_{0\le s\le 1, 0\le t-s\le 1}|B_{s+(t-s)}-B_s|\geq \frac{Cc_0}{ \eta}\sqrt{1-y\eta^2} \right).
\end{eqnarray*}
Following \cite[Theorem 12.1.c]{lin2011probability}, we can get
\begin{eqnarray*}
\sum_{k=0}^{m^\alpha-1}\PP\left(\sup_{s\in[m-k-1,m-k]}(\tilde{\mathcal{H}}_{m-k}-\tilde{\mathcal{H}}_{s})  \geq \frac{c_0}{2}\sqrt{1-y\eta^2} \right)
&\le& cm^\alpha e^{-\frac{Cc_0^2}{\eta^2}(1-y\eta^2)}.
\end{eqnarray*}
Hence we have
\begin{eqnarray}\label{e:sum8}
&&\PP\left(\frac{\tilde{\mathcal{H}}_m-\tilde{\mathcal{H}}_{T_1}}{\sqrt{1-y\eta^2}} \geq c_0, T_1\in [m-m^\alpha,m+m^\alpha] \right)\\
&\le& (1+m^\alpha)\big(e^{-\frac{C(y'-m^\alpha)^2}{m^\alpha}}+e^{-\frac{c_0^2(m-y)}{8y'}}\big)+cm^\alpha e^{-C \eta^{-2}c_0^2(1-y\eta^2)}.\nonumber
\end{eqnarray}

Combining (\ref{e:sum4}-\ref{e:sum8}), we obtain
\begin{eqnarray}
\PP\left(\mathcal{H}_\eta/\sqrt{\Y_\eta} \geq x\right)
&\le& 1-\Phi(\sqrt{1-y\eta^2}(x-c_0))+2e^{-Cy^2\eta^2}+2e^{-Cm^{2\alpha-1}}\nonumber\\
&&+(1+m^\alpha)\big(e^{-\frac{C(y'-m^\alpha)^2}{m^\alpha}}+e^{-\frac{c_0^2(m-y)}{8y'}}\big)+cm^\alpha e^{-C \eta^{-2}c_0^2(1-y\eta^2)}.
\end{eqnarray}
By the following well known estimate of normal distribution (cf. \citet[(4.1)]{Fan2019Self})
\begin{align}\label{bound of normal}
\frac{1}{\sqrt{2\pi}(1+x)}e^{-\frac{x^2}{2}}\le 1-\Phi(x)\le\frac{1}{\sqrt{\pi}(1+x)}e^{-\frac{x^2}{2}},~~x\ge0,
\end{align}
we have
\begin{eqnarray*}
&&\PP\left(\mathcal{H}_\eta/\sqrt{\Y_\eta} \geq x\right) \Big/(1-\Phi(x))\\
&\le& \frac{1-\Phi(\sqrt{1-y\eta^2}(x-c_0))}{1-\Phi(x)}+\sqrt{2\pi}(1+x)\left[2e^{-Cy^2\eta^2+\frac{x^2}{2}}+2e^{-Cm^{2\alpha-1}+\frac{x^2}{2}}\right.\\
&&\left. +(1+m^\alpha)\big(e^{-\frac{C(y'-m^\alpha)^2}{m^\alpha}+\frac{x^2}{2}}+e^{-\frac{c_0^2(m-y)}{8y'}+\frac{x^2}{2}}\big)+cm^\alpha e^{-C \eta^{-2}c_0^2(1-y\eta^2)+\frac {x^2}{2}}\right].
\end{eqnarray*}
For the normal distribution part, by \eqref{bound of normal} again, we have
\begin{eqnarray}\label{e:compareN}
&&\frac{1-\Phi(\sqrt{1-y\eta^2}(x-c_0))}{1-\Phi(x)}\\
&=&\frac{1-\Phi(\sqrt{1-y\eta^2}(x-c_0))}{1-\Phi(x)}1_{\{x\ge1\}}+\frac{1-\Phi(\sqrt{1-y\eta^2}(x-c_0))}{1-\Phi(x)}1_{\{0\le x<1\}}\nonumber\\
&\le&\frac{\sqrt2(1+x)}{1+\sqrt{1-y\eta^2}(x-c_0)}e^{\frac{1}{2}x^2-\frac{1}{2}(1-y\eta^2)(x-c_0)^2}1_{\{x\ge1\}}\nonumber\\
&&+\big(1+\int_{\sqrt{1-y\eta^2}(x-c_0)}^xe^{-\frac{t^2}{2}}\dif t\big/\int_x^\infty e^{-\frac{t^2}{2}}\dif t\big)1_{\{0\le x<1\}}\nonumber\\
&\le&\frac{\sqrt2(1+x)}{1+\sqrt{1-y\eta^2}(x-c_0)}e^{-\frac{1}{2}c_0^2+xc_0+\frac12(x-c_0)^2y\eta^2}1_{\{x\ge1\}}\nonumber\\
&&+\big[1+(1+x)\big(x-\sqrt{1-y\eta^2}(x-c_0)\big)e^{-\frac{1}{2}c_0^2+xc_0+\frac12(x-c_0)^2y\eta^2}\big]1_{\{0\le x<1\}}\nonumber\\
&\le&e^{-\frac{1}{2}c_0^2+xc_0+C(x-c_0)^2y\eta^2}1_{\{x\ge1\}}+\big[1+C\big(x-\sqrt{1-y\eta^2}(x-c_0)\big)\big]1_{\{0\le x<1\}}.\nonumber
\end{eqnarray}
Thus,
\begin{eqnarray*}
&&\PP\left(\mathcal{H}_\eta/\sqrt{\Y_\eta} \geq x\right) \Big/(1-\Phi(x))\\
&\le&e^{-\frac{1}{2}c_0^2+xc_0+C(x-c_0)^2y\eta^2}1_{\{x\ge1\}}+\big[1+C\big(x-\sqrt{1-y\eta^2}(x-c_0)\big)\big]1_{\{0\le x<1\}}\\
&&+\sqrt{2\pi}(1+x)\left[2e^{-Cy^2\eta^2+\frac{x^2}{2}}+2e^{-Cm^{2\alpha-1}+\frac{x^2}{2}}\right.\\
&&\left.+(1+m^\alpha)\big(e^{-\frac{C(y'-m^\alpha)^2}{m^\alpha}+\frac{x^2}{2}}+e^{-\frac{c_0^2(m-y)}{8y'}+\frac{x^2}{2}}\big)+cm^\alpha e^{-C \eta^{-2}c_0^2(1-y\eta^2)+\frac {x^2}{2}}\right].
\end{eqnarray*}
To guarantee the limit of the first two terms is $1$ and the last term is $0$ as $\eta\to0$, i.e. $m\to\infty$. We need $y^2\eta^2\to\infty$, $2Cy^2\eta^2>x^2$, $\eta^2x^2y\to0$ and $x=o(c_0^{-1})$. Choosing $y=\eta^{-\frac43}$, $c_0=\eta^{\frac13}$, $y'=\eta^{-\frac23}$ and $\alpha=2/3$, one has
\begin{eqnarray*}
&&\PP\left( \frac{\mathcal{H}_\eta}{\sqrt{\Y_\eta}} \geq x\right) \Big/(1-\Phi(x))\\
&\le& e^{-c(\eta^{\frac23}-x\eta^{\frac13}-x^2\eta^{\frac23})}1_{\{x\ge1\}}+\big[1+C\big(x-(1-\eta^{\frac23})^{-\frac12}(x-\eta^{\frac13})\big)\big]1_{\{0<x<1\}}+e^{C(x^2-\eta^{-\frac23})}\\
&\le& 1+C(x\eta^{\frac13}+\eta^{\frac13})
\end{eqnarray*}
converges to $1$ uniformly for $\eta^{\frac13}\le x=o(\eta^{-1/3})$ as $\eta$ tends to $0$.

For the lower bound of $\PP\left( \frac{\mathcal{H}_\eta}{\sqrt{\Y_\eta }} \geq x\right) \Big/(1-\Phi(x))$, we have
\begin{eqnarray*}
\PP\left( \frac{\mathcal{H}_\eta}{\sqrt{\Y_\eta }} \geq x\right)
&\ge& \PP\left( \frac{\mathcal{H}_\eta}{\sqrt{\Y_\eta }} \geq x, \eta^{-2}\left|1 - \Y_\eta\right|\le y\right)\\
&\ge& \PP\left( \frac{\mathcal{H}_\eta}{\sqrt{1+\eta^{2}y}} \geq x, \eta^{-2}\left|1 - \Y_\eta\right|\le y\right)\\
&=&\PP\left( \frac{\mathcal{H}_\eta}{\sqrt{1+\eta^{2}y}} \geq x\right)-\PP\left( \frac{\mathcal{H}_\eta}{\sqrt{1+\eta^{2}y}} \geq x, \eta^{-2}\left|1 - \Y_\eta\right|> y\right)\\
&\ge& \PP\left( \frac{\tilde{\mathcal{H}}_{T_1}}{\sqrt{1+\eta^{2}y}} \geq x+c_0\right)-\PP\left( \frac{\tilde{\mathcal{H}}_{T_1}-\tilde{\mathcal{H}}_m}{\sqrt{1+\eta^{2}y}} \geq c_0\right)-\PP\left( \eta^{-2}\left|1 - \Y_\eta\right|> y\right).
\end{eqnarray*}
Similar with the estimate of the upper bound, (\ref{bound of normal}), (\ref{e:6}) and Lemma \ref{EY-Y} imply
\begin{eqnarray*}
\PP\left( \frac{\mathcal{H}_\eta}{\sqrt{\Y_\eta}} \geq x\right) \Big/(1-\Phi(x))
\ge 1-C(x\eta^{\frac13}+\eta^{\frac13})
\end{eqnarray*}
converges to $1$ uniformly for $\eta^{\frac13}\le x=o(\eta^{-1/3})$ as $\eta$ tends to $0$. Hence, we have
\begin{eqnarray*}
\PP\left( \frac{\mathcal{H}_\eta}{\sqrt{ \Y_\eta }} \geq x\right) \Big/(1-\Phi(x))
&=& 1+O(x\eta^{\frac13}+\eta^{\frac13})
\end{eqnarray*}
uniformly for $\eta^{\frac13}\le x=o(\eta^{-1/3})$ as $\eta$ vanishes.

\end{proof}

\subsection{Proof of Theorem \ref{thm-self-normal}}
\Bp[Proof of  Theorem \ref{thm-self-normal}]
We have proved the following decomposition,
$$
\eta^{-\frac12}\big( \Pi_\eta(h) - \pi(h) \big)= \mathcal{R}_{\eta}+\mathcal{H}_{\eta}.
$$
Noting that, for any $x > 0$ and $0<y<x$, we have
\begin{align}
\PP(\W_{\eta} \geq x) & =\PP\left( \frac{\mathcal{R}_{\eta}+\mathcal{H}_{\eta}}{\sqrt{ \Y_\eta}} \geq x\right)\leq \PP\left( \frac{\mathcal{H}_\eta}{\sqrt{ \Y_\eta}} \geq x-y\right) + \PP\left( \frac{\mathcal{R}_\eta}{\sqrt{ \Y_\eta}} \geq y\right).\label{eq-PW upper bound}
\end{align}
For the first term, Lemma \ref{Fan_Coro2.1} yields that
$$
\PP\left( \frac{\mathcal{H}_\eta}{\sqrt{ \Y_\eta}} \geq x-y\right) \Big/(1-\Phi(x-y))  =1+O((x-y)\eta^{\frac13}+\eta^{\frac13})
$$
uniformly for $\eta^{\frac13}\le x-y =o(\eta^{-\frac 13})$ as $\eta$ tends to zero. We take $\eta^{\frac13}<x=o(\eta^{-\alpha})$ and $y=o(1)$ such that $\alpha\le1/3$, $xy\to0$ and $x-y\ge\eta^{\frac13}$, here $y$ will be chosen later. Similar with the calculation of \eqref{e:compareN}, (\ref{bound of normal}) yields
\begin{align*}
\frac{1-\Phi(x-y)}{1-\Phi(x)}= 1+O(xy+y).
\end{align*}
Hence,
\begin{align}\label{mdp_1}
\frac{\PP\left( \mathcal{H}_\eta/\sqrt {\Y_\eta} \geq x-y\right)} {1-\Phi(x)}
&=\frac{\PP\left(\mathcal{H}_\eta/\sqrt {\Y_\eta} \geq x-y\right)}{1-\Phi(x-y)}\frac{1-\Phi(x-y)}{1-\Phi(x)  }\nonumber\\
&= 1+O(x\eta^{\frac13}+\eta^{\frac13}+xy+y)
\end{align}
as $\eta$ vanishes.

For the second term of  (\ref{eq-PW upper bound}), we have

\begin{eqnarray*}
\PP\left( \frac{\mathcal{R}_\eta}{\sqrt{\Y_\eta}} \geq y\right)
&\le&\PP\left(\Y_\eta<\E\Y_\eta-y\right)+\PP\left(\frac{\mathcal{R}_\eta}{\sqrt{\Y_\eta}} \geq y, \Y_\eta \ge \E\Y_\eta-y\right)\nonumber\\
&\le&\PP\left(\E\Y_\eta-\Y_\eta> y\right)+\PP\left( \mathcal{R}_\eta\ge y\sqrt{\E\Y_\eta-y}\right).
\end{eqnarray*}
For the first probability, Lemma \ref{EY-Y} yields that
\begin{align*}
\PP(\E\Y_\eta-\Y_\eta>y)\le e^{-Cy^2\eta^{-2}}.
\end{align*}
For the second probability, following the stationary of $\theta_k$ and Lemma \ref{lem-R}, one has
\begin{align*}
\PP\left( \mathcal{R}_\eta\ge y\sqrt{\E\Y_\eta-y}\right)&=\PP\left( \mathcal{R}_{\eta}\ge y\sqrt{\pi(|\sigma^{\mathrm{T}}\nabla \vp|^2)-y}\right)
\le C e^{-c\eta^{-2\bar\gamma}y^{\frac 23}},
\end{align*}
as $y\ge c\max\{\eta^{\frac 32-6\bar\gamma},\eta^{\frac32\bar\gamma},\eta^{\frac 12}\}=c\eta^{\frac 32-6\bar\gamma}$ where $\frac15\le\bar\gamma<\frac14$. Hence, we have
\begin{eqnarray*}
\PP\left( \frac{\mathcal{R}_\eta}{\sqrt{\Y_\eta}} \geq y\right)
&\le&C\left(e^{-c\eta^{-2}y^2}
+e^{-c\eta^{-2\bar\gamma}y^{\frac 23}}\right).
\end{eqnarray*}
This, together with (\ref{bound of normal}), implies
\begin{align*}
\PP\left( \frac{\mathcal{R}_{\eta}}{\sqrt{\Y_\eta}} \geq y\right) \Big/(1-\Phi(x))
\le&C(1+x)e^{\frac12 x^2}\PP\left( \frac{\mathcal{R}_\eta}{\sqrt{\Y_\eta}} \geq y\right).
\end{align*}
It converges to $0$ as $\eta\to0$ uniformly for
\begin{align*}
\eta^{\frac13}\le x=o\left(\min\{\eta^{-1}y,\eta^{-\bar\gamma}y^{\frac13}
\}\right).
\end{align*}
Since Lemma \ref{Fan_Coro2.1} holds uniformly as $\eta^{\frac13}+y\le x=o(\eta^{-\alpha})$, we need to choose $\alpha$, $y$ and $\bar\gamma$ such that
$$\min\{
\eta^{-1}y, \eta^{-\bar\gamma}y^{\frac13}\}\ge \eta^{-\alpha}.$$
By taking  $\alpha=1/6$, $y=c\eta^{\frac16}$ and $\bar\gamma=2/9$, we can get
\begin{align}\label{mdp_2}
\PP\left( \frac{\mathcal{R}_{\eta}}{\sqrt{\Y_\eta}} \geq y\right) \Big/(1-\Phi(x))
\le&C(1+x)\exp\{c(x^2-\eta^{-\frac13})\}\to0
\end{align}
uniformly for $c\eta^{\frac16}\le x=o(\eta^{-\frac 16})$ as $\eta$ vanishes.

Following  (\ref{eq-PW upper bound}), (\ref{mdp_1}) and (\ref{mdp_2}), we have
\begin{align}\label{eq-PW upper bound 1}
&\PP\left( \frac{\mathcal{R}_\eta+\mathcal{H}_\eta}{\sqrt \Y_\eta} \geq x\right)\Big/(1-\Phi(x))\le 1+C(x\eta^{1/6}+\eta^{1/6}).
\end{align}
uniformly for $c\eta^{\frac 16}\le x=o(\eta^{-\frac 16})$ as $\eta$ tends to zero.

On the other hand,
\begin{align*}
\PP\left( \frac{\mathcal{R}_\eta+\mathcal{H}_\eta}{\sqrt{\Y_\eta}} \geq x\right)
\ge \PP\left( \frac{\mathcal{H}_\eta}{\sqrt{\Y_\eta}} \geq x+y\right) - \PP\left( \frac{-\mathcal{R}_\eta}{\sqrt{\Y_\eta}} \geq y\right).
\end{align*}
Similar as the proof of  (\ref{eq-PW upper bound 1}), Lemmas \ref{lem-R} and \ref{Fan_Coro2.1} yield that
\begin{align*}
&\PP\left( \frac{\mathcal{R}_\eta+\mathcal{H}_\eta}{\sqrt{\Y_\eta}} \geq x\right)\Big/(1-\Phi(x))\geq 1-C(x\eta^{\frac 16}+\eta^{\frac16}),
\end{align*}
uniformly for $c\eta^{\frac16}\le x=o(\eta^{-\frac 16})$ as $\eta$ tends to zero. Combining the last inequality with (\ref{eq-PW upper bound 1}), we deduce that
\begin{align*}
&\PP\left( \W_\eta \geq x\right)\Big/(1-\Phi(x))=1+ O\left(x\eta^{1/6}+\eta^{1/6}\right),
\end{align*}
uniformly for $c\eta^{\frac16}\le x=o(\eta^{-\frac 16})$ as $\eta$ tends to zero.
\Ep

\begin{appendix}

\section{Proofs of Lemmas in Section \ref{s:2}}
\begin{proof}[Proof of Lemma \ref{l:Erg2}]
We first give the proof of the ergodicity of $(X_t)_{t\ge0}$. Following \citet[Theorem 2.1]{Tweedie1996}, it is easy to verify the irreducibility of $(X_t)_{t\ge0}$. For the Lyapunov function $V(x)={|x|^2+1}$, following (\ref{assu4}) and (\ref{generator}), we have
\begin{eqnarray*}
\mathcal{A}V(x)=\Ll g(x),2x\Rr+\|\sigma\|^2\le -K_1|x|^2+C\le-\frac{K_1}{2}V(x)+(C+\frac{K_1}{2})1_{\{|x|\le\sqrt{2C/K_1+1}\}}.
\end{eqnarray*}
By \citet[Theorem 6.1]{MeTw93}, $(X_t)_{t\ge0}$ is exponential ergodic with invariant measure $\pi$ satisfying
\begin{align}\label{ergodic1}
\left|\E[h(X^x_t)-\pi(h)]\right|\le CV(x) e^{-ct}.
\end{align}

Then we consider the ergodicity of $(\theta_k )_{k\ge0}$. Denote its transition probability by $P(x,\dif y)$ for $x,y\in\R^d$. For any open set $A\in\R^d$ and initial value $x$, since $\xi_1$ is a normal random vector, we have
\begin{align*}
P(x,A)=\PP(x+\eta g(x)+\sqrt\eta\sigma\xi_1\in A)>0.
\end{align*}
Suppose $P^k(x,A)>0$ for some integer $k> 1$, then we have
\begin{align*}
P^{k+1}(x,A)=\int_{\R^d}P(x,y)P^k(y,A) \dif y>0.
\end{align*}
The induction yields that $(\theta_k )_{k\ge0}$ is irreducible. Following (\ref{assu3}),(\ref{assu4}) and (\ref{Lang}), one has
\begin{eqnarray}\label{e:lya}
\E_k[V(\theta_{k+1})]&=&\E_k[|\theta_k +\eta g(\theta_k )+\sqrt\eta\sigma\xi_{k+1}|^2]+1\nonumber\\
&=&|\theta_k|^2+|\eta g(\theta_k )|^2+\eta \|\sigma\|^2+2\Ll\theta_k ,\eta g(\theta_k )\Rr+1\nonumber\\
&\le&(1-K_1\eta+2L^2\eta^2)|\theta_k |^2+2|g(0)|^2\eta^2+\eta\|\sigma\|^2+2C\eta+1\nonumber\\
&\le&(1-\frac{1}{2}K_1\eta+2L^2\eta^2)V(\theta_k )+b1_{D}(\theta_k ).
\end{eqnarray}
Here $b=\frac{1}{2}K_1\eta-2L^2\eta^2+2|g(0)|^2\eta^2+\eta\|\sigma\|^2+2C\eta$ and set $D=\{|x|\le \frac{2b}{\eta K_1}\}$. There exists $\eta_0$ such that for $\eta\le\eta_0$, $1-\frac{1}{2}K_1\eta+2L^2\eta^2<1$. By \citet[(29)]{Tweedie1996}, we deduce that $\theta_k $ is ergodic when $\eta\le\eta_0$, that is
\begin{align}\label{ergodic2}
\left|\E[h(\theta_k )-\pi_\eta(h)]\right|\le CV(\theta_0)e^{-ck}.
\end{align}
\end{proof}

Moreover, for the function $\tilde V(x)=|x|^4+1$, similarly with the calculation of (\ref{e:lya}), we can get
\begin{eqnarray*}
\E_k[\tilde V(\theta_{k+1})]&=&\E_k[|\theta_{k+1}|^4+1]\\
&\le&(1-2K_1\eta+C_1\eta^2)|\theta_k|^4+C_2\eta|\theta_k|^2+C_3\eta^2+1,
\end{eqnarray*}
where $C_1,C_2,C_3$ depend on $\sigma$, $K_1$, $L$ and $C$ in (\ref{assu3}), (\ref{assu4}). Then we have
\begin{eqnarray}  \label{e:BenDan-0}
\E_k[\tilde V(\theta_{k+1})]&\le&(1-K_1\eta+C_1\eta^2)\tilde V(\theta_k)+\tilde b1_{\tilde D}({\theta_k}),
\end{eqnarray}
where $\tilde b=\frac{C_2^2}{4K_1}\eta+K_1\eta+C_3\eta^2$, $\tilde D=\{|x|^2\le(\frac{C_3-C_1}{K_1}\eta+1+(\frac{C_2}{2K_1})^2)^{\frac12}+\frac{C_2}{2K_1}\}$. For small enough $\eta$ such that $1-K_1\eta+C_1\eta^2<1$, let $\theta_0$ take the ergodic measure $\pi_\eta$, then $(\theta_k)_{k \ge 0}$ is stationary and \eqref{e:BenDan-0} implies
\begin{eqnarray*}  \label{e:BenDan-1}
\pi_\eta(\tilde V)&\le&(1-K_1\eta+C_1\eta^2)\pi_\eta(\tilde V)+\tilde b,
\end{eqnarray*}
i.e.
\begin{eqnarray}  \label{e:BenDan-1}
\pi_\eta(\tilde V)&\le& \frac{\tilde b}{K_1\eta-C_1\eta^2}.
\end{eqnarray}

Notice that for any $k=0,...,m$ and positive number $\gamma$, we have
\begin{eqnarray*}
&&\E_k\left[e^{\gamma|\theta_{k+1}|^2}\right]\\
&=&\E_k\left[\exp\{\gamma|\theta_k|^2+\gamma|\eta g(\theta_k )|^2+\gamma\eta |\sigma\xi_{k+1}|^2+2\gamma\Ll\theta_k ,\eta g(\theta_k )\Rr+2\sqrt\eta\gamma\Ll\sigma^{\rm {T}}(\theta_k+\eta g(\theta_k)) ,\xi_{k+1}\Rr\}\right]\\
&=& e^{\gamma|\theta_k|^2+\gamma|\eta g(\theta_k )|^2+2\gamma\Ll\theta_k ,\eta g(\theta_k )\Rr}\E_k\left[\exp\{\gamma\eta |\sigma\xi_{k+1}|^2+2\sqrt\eta\gamma\Ll\sigma^{\rm {T}}(\theta_k+\eta g(\theta_k)) ,\xi_{k+1}\Rr\}\right].
\end{eqnarray*}
A straight calculation to the conditional expectation with respect to the Gaussian random variable $\xi_{k+1}$  yields
\begin{eqnarray*}
\E_k\left[\exp\{\gamma\eta |\sigma\xi_{k+1}|^2+2\sqrt\eta\gamma\Ll\sigma^{\rm {T}}(\theta_k+\eta g(\theta_k)) ,\xi_{k+1}\Rr\}\right]
\le2\exp\left\{4\eta\gamma^2\|\sigma\|^2|\theta_k+\eta g(\theta_k)|^2\right\},
\end{eqnarray*}
here $\gamma$ is chosen small enough such that $\gamma\|\sigma\|^2\le1/4$. This estimate, together with (\ref{assu3}) and (\ref{assu4}), implies
\begin{eqnarray*}
\E_k\left[e^{\gamma|\theta_{k+1}|^2}\right]&\le& 2\exp\left\{(1-K_1\eta+4\eta\gamma+C_1\eta^2)\gamma|\theta_k|^2+3C\gamma\eta\right\}\\
&\le& (1-K_1\eta+4\eta\gamma+C_1\eta^2) e^{\gamma|\theta_k|^2}+\bar b,
\end{eqnarray*}
with $\eta$ and $\gamma$ are small enough such that $1-K_1\eta+4\eta\gamma+C_1\eta^2<1$ and $\bar b$ is big enough such that the second inequality holds. Let $\theta_0$ take the ergodic measure $\pi_\eta$, then we have
\begin{eqnarray}  \label{e:e-lyapunov}
\pi_\eta(e^{\gamma|\cdot|^2})&\le& \frac{\bar b}{K_1\eta-4\eta\gamma-C_1\eta^2}.
\end{eqnarray}

\section{The proof of lemmas in section 4 }
\begin{proof}[Proof of Lemma \ref{C2}]
For the first inequality, by using H\"{o}lder's inequality, we can get
\begin{eqnarray*}
&&\E_k\exp\left\{\Ll \Psi_1(\theta_k),\sigma\xi_{k+1}   \Rr+\Psi_2(\theta_k,\xi_{k+1})\right\}\\
&\le&\left(\E_k\exp\{2\Ll \Psi_1(\theta_k),\sigma\xi_{k+1} \Rr-2|\sigma^{\mathrm{T}}\Psi_1(\theta_k)|^2 \}\right)^{\frac 12}
\left(\E_k\exp\{2\Psi_2(\theta_k,\xi_{k+1})+2|\sigma^{\mathrm{T}}\Psi_1(\theta_k)|^2 \}\right)^{\frac 12}.
\end{eqnarray*}
Since $\xi_{k+1}$ is gaussian distributed and independent of $\theta_k$, a straightforward calculation gives
\begin{align*}
\left(\E_k\exp\{2\Ll \Psi_1(\theta_k),\sigma\xi_{k+1} \Rr-2|\sigma^{\mathrm{T}}\Psi_1(\theta_k)|^2   \}\right)^{\frac 12}=1.
\end{align*}
Hence, we have
\begin{eqnarray*}
\E_k\exp\left\{\Ll \Psi_1(\theta_k),\sigma\xi_{k+1}   \Rr+\Psi_2(\theta_k,\xi_{k+1})\right\}
\le \left(\E_k \exp\left\{2|\Psi_1(\theta_k)|^2\|\sigma\|^2+2\Psi_2(\theta_k,\xi_{k+1})\right\}\right)^{\frac 12}.
\end{eqnarray*}

For the second inequality of Lemma \ref{C2}, by the same way we have
\begin{eqnarray*}
&&\E \exp\left\{\sum_{k=0}^{m-1}\left(\Ll \Psi_1(\theta_k),\sigma\xi_{k+1}   \Rr+\Psi_2(\theta_k,\xi_{k+1})\right)\right\}\\
&\le&\left(\E \exp\left\{\sum_{k=0}^{m-1}2\left(\Ll \Psi_1(\theta_k),\sigma\xi_{k+1}   \Rr-|\sigma^{\mathrm{T}}\Psi_1(\theta_k)|^2\right)\right\}\right)^{\frac 12}\\
&&\times\left(\E \exp\left\{\sum_{k=0}^{m-1}2\left(|\sigma^{\mathrm{T}}\Psi_1(\theta_k)|^2+\Psi_2(\theta_k,\xi_{k+1})\right)\right\}\right)^{\frac12}\\
&=& \left(\E \exp\left\{\sum_{k=0}^{m-1}2\left(|\Psi_1(\theta_k)|^2\|\sigma\|^2+\Psi_2(\theta_k,\xi_{k+1})\right)\right\}\right)^{\frac 12},
\end{eqnarray*}
where the following relation is obtained by a standard conditional argument:
\begin{eqnarray*}
&&\E\exp\left\{\sum_{k=0}^{m-1}2\left(\Ll \Psi_1(\theta_k),\sigma\xi_{k+1}   \Rr-|\sigma^{\mathrm{T}}\Psi_1(\theta_k)|^2\right)\right\}\\
&=&\E\left[\exp\left\{\sum_{k=0}^{m-2}2\left(\Ll \Psi_1(\theta_k),\sigma\xi_{k+1}   \Rr-|\sigma^{\mathrm{T}}\Psi_1(\theta_k)|^2\right)\right\}\E_{m-1}\left[e^{2\Ll \Psi_1(\theta_{m-1}),\sigma\xi_{m}   \Rr-2|\sigma^{\mathrm{T}}\Psi_1(\theta_{m-1})|^2 }\right]\right]\\
&=&\E \exp\left\{\sum_{k=0}^{m-2}2\left(\Ll \Psi_1(\theta_k),\sigma\xi_{k+1}   \Rr-|\sigma^{\mathrm{T}}\Psi_1(\theta_k)|^2\right)\right\}=...=1.
\end{eqnarray*}
A similar calculation gives the third inequality.

\end{proof}

\begin{proof}[Proof of Lemma \ref{lem-g}]
Since $\theta_{k+1} =\theta_k +\eta g(\theta_k )+\sqrt\eta\sigma\xi_{k+1}$, it is easy to calculate that
$$|\theta_{k+1} |^2-|\theta_k |^2=\eta^2 |g(\theta_k )|^2+\eta|\sigma\xi_{k+1}|^2+2\Ll\eta\theta_k  ,g(\theta_k )\Rr+2\Ll\sqrt\eta\theta_k+\eta^{\frac32}g(\theta_k ) ,\sigma\xi_{k+1}\Rr.$$
Summing these equalities from $k=0$ to $k=m-1$, we obtain
\begin{align}\label{2}
|\theta_{m} |^2-|\theta_0|^2=\sum_{k=0}^{m-1}\left[\eta^2 |g(\theta_k )|^2+\eta|\sigma\xi_{k+1}|^2+2\Ll\eta\theta_k  ,g(\theta_k )\Rr+2\Ll\sqrt\eta\theta_k+\eta^{\frac32}g(\theta_k ) ,\sigma\xi_{k+1}\Rr\right].
\end{align}
For $\gamma>0$, (\ref{assu4}) and (\ref{2}) imply
\begin{eqnarray*}
&&\E_0\exp\left\{\sum_{k=0}^{m-1}\frac{K_1}{2}\gamma\eta|\theta_k |^2\right\}\le \E_0\exp\left\{-\sum_{k=0}^{m-1}\gamma\Ll\eta\theta_k , g(\theta_k )\Rr\right\}e^{C\gamma\eta^{-1}}\\
&\le&\E_0\exp\left\{\frac{\gamma|\theta_0|^2}{2}+\frac{\gamma}{2}\sum_{k=0}^{m-1}\left[\eta^2 |g(\theta_k )|^2+\eta|\sigma\xi_{k+1}|^2+2\Ll\sqrt\eta\theta_k+\eta^{\frac32}g(\theta_k ) ,\sigma\xi_{k+1}\Rr\right]\right\}e^{C\gamma\eta^{-1}}\\
&=&\E_0 \exp\left\{\sum_{k=0}^{m-1}\left[\frac{\gamma}{2}\left(\eta^2 |g(\theta_k )|^2+\eta|\sigma\xi_{k+1}|^2\right)+\gamma\Ll\sqrt\eta\theta_k+\eta^{\frac32}g(\theta_k ) ,\sigma\xi_{k+1}\Rr\right]\right\}e^{\frac{\gamma|\theta_0|^2}{2}+C\gamma\eta^{-1}}.
\end{eqnarray*}
By Lemma \ref{C2} with
$\Psi_1(\theta_k)= \gamma(\sqrt\eta\theta_k+\eta^{\frac32}g(\theta_k ))$ and
$\Psi_2(\theta_k,\xi_{k+1})= \frac{\gamma}{2}(\eta^2 |g(\theta_k )|^2+\eta|\sigma\xi_{k+1}|^2)$ therein, we have
\begin{eqnarray*}
&&\E_0 \exp\left\{\sum_{k=0}^{m-1}\left[\frac{\gamma}{2}\left(\eta^2 |g(\theta_k )|^2+\eta|\sigma\xi_{k+1}|^2\right)+\gamma\Ll\sqrt\eta\theta_k+\eta^{\frac32}g(\theta_k ) ,\sigma\xi_{k+1}\Rr\right]\right\}\\
&\le&\left(\E_0\exp\left\{\sum_{k=0}^{m-1}\left(\gamma(\eta^2|g(\theta_k)|^2+\eta|\sigma\xi_{k+1}|^2)
+2\gamma^2|\sqrt\eta\theta_k+\eta^{\frac32}g(\theta_k )|^2\|\sigma\|^2
\right)\right\}\right)^{\frac 12}\\
&\le&\left(\E_0\exp\left\{\sum_{k=0}^{m-1}2\gamma\eta|\sigma\xi_{k+1}|^2\right\}\right)^{\frac14}
\left(\E_0\exp\left\{\sum_{k=0}^{m-1}(2\gamma\eta^2|g(\theta_k )|^2+4\gamma^2|\sqrt\eta\theta_k+\eta^{\frac32}g(\theta_k )|^2\|\sigma\|^2)\right\}\right)^{\frac14}.
\end{eqnarray*}
For the first expectation, we take some $\gamma_0'$ and $\eta_0'$ such that $1-4\gamma_0'\eta_0'\|\sigma\|^2>0$. Then for any $\gamma<\gamma_0'$ and $\eta<\eta_0'$, we have
\begin{eqnarray*}
\E_0\exp\left\{\sum_{k=0}^{m-1}2\gamma\eta|\sigma\xi_{k+1}|^2\right\}
&\le& \E_0\exp\left\{\sum_{k=0}^{m-1}2\gamma\eta\|\sigma\|^2|\xi_{k+1}|^2\right\}\\
&=&\left(\int_{-\infty}^\infty\frac{1}{\sqrt{2\pi}}\exp\{2\gamma\eta\|\sigma\|^2 x^2-\frac{1}{2}x^2\}\dif x\right)^{md}\\
&=&(1-4\gamma\eta\|\sigma\|^2)^{-\frac {md}2}.
\end{eqnarray*}
For the second expectation, by (\ref{assu3}), we can choose some $\gamma_0''$ and $\eta_0''$ such that as $\gamma<\gamma_0''$ and $\eta<\eta_0''$
\begin{align}\label{eq-lem4.2}
2\gamma\eta^2|g(\theta_k )|^2+4\gamma^2|\sqrt\eta\theta_k+\eta^{\frac 32}g(\theta_k)|^2\|\sigma\|^2\le \frac{K_1}{2}\gamma\eta|\theta_k|^2+C\eta\gamma,
\end{align}
which leads to
\begin{eqnarray*}
\E_0\exp\left\{\sum_{k=0}^{m-1}(2\gamma\eta^2|g(\theta_k )|^2+4\gamma^2|\sqrt\eta\theta_k+\eta^{\frac32}g(\theta_k )|^2\|\sigma\|^2)\right\}\le e^{C \gamma \eta^{-1}}\E_0 \exp\left\{\sum_{k=0}^{m-1}\frac{K_1}{2}\gamma\eta|\theta_k|^2\right\}.
\end{eqnarray*}
Hence, for $\gamma<\gamma_0=\gamma_0'\wedge\gamma_0''$ and $\eta<\eta_0=\eta_0'\wedge\eta_0''$, we have
\begin{eqnarray*}
\E_0\exp\left\{\sum_{k=0}^{m-1}\frac{K_1}{2}\gamma\eta|\theta_k |^2\right\}\le e^{\frac{\gamma|\theta_0|^2}{2}+C\gamma\eta^{-1}} (1-4\gamma\eta\|\sigma\|^2)^{-\frac {md} 8}\left(\E_0\exp\left\{\sum_{k=0}^{m-1}\frac{K_1}{2}\gamma\eta|\theta_k |^2\right\}\right)^{\frac 14},
\end{eqnarray*}
i.e.,
\begin{align*}
\left(\E_0\exp\left\{\sum_{k=0}^{m-1}\frac{K_1}{2}\gamma\eta|\theta_k |^2\right\}\right)^{\frac 34}\le e^{\frac{\gamma|\theta_0|^2}{2}+C\gamma\eta^{-1}} (1-4\gamma\eta\|\sigma\|^2)^{-\frac {md} 8}.
\end{align*}
Then we have
\begin{eqnarray}
\E_0\exp\left\{\sum_{k=0}^{m-1}\frac{K_1}{2}\gamma\eta|\theta_k |^2\right\}&\le&  e^{\frac{2\gamma|\theta_0|^2}{3}+C\gamma\eta^{-1}}\left((1-4\gamma\eta\|\sigma\|^2)^{-\frac{1}{4\gamma\eta\|\sigma\|^2}}\right)^{\frac{2\gamma\|\sigma\|^2 d} {3\eta}}\nonumber\\
&\le& e^{\frac{2\gamma|\theta_0|^2}{3}+C \gamma\eta^{-1}}.  \label{e:ConExp}
\end{eqnarray}
This, together with \eqref{assu3}, implies
\begin{eqnarray*}
\E_0 \exp\left\{\sum_{k=0}^{m-1}\frac{K_1}{4L^2}\gamma\eta|g(\theta_k)|^2\right\}
&\le&\E_0 \exp\left\{\sum_{k=0}^{m-1}\frac{K_1}{2}\gamma\eta|\theta_k |^2\right\}  e^{\frac{K_1 |g(0)|^2}{2L^2} \gamma\eta^{-1}}\\
&\le& C e^{c (\eta^{-1}+|\theta_0|^2)}.
\end{eqnarray*}
Writing $\tl \gamma=\frac{K_1}{4L^2}\gamma$ and replacing the $\gamma$ in \eqref{R_1,22-0} by $\tl \gamma$, we immediately finish the proof of \eqref{R_1,22-0}.

For \eqref{R_1,22} with $\theta_0\sim\pi_\eta$, \eqref{R_1,22-0} and \eqref{e:e-lyapunov} yield
\begin{eqnarray*}
\E\exp\left\{\gamma \eta\sum_{k=0}^{m-1}|g(\theta_k)|^2\right\} &=&\E\left[\E_0 \exp\left\{\gamma \eta\sum_{k=0}^{m-1}|g(\theta_k)|^2\right\}\right]\le C e^{c \eta^{-1}}
\end{eqnarray*}

The inequalities \eqref{e:R-1-22-0} and \eqref{e:R-1-22} immediately follow by Chebyshev's inequality.
\end{proof}

\Bp[Proof of Lemma \ref{C5}]
It is easy to see that
\begin{eqnarray*}
\E \exp\left\{\frac{1}{\sqrt m}\sum_{n=0}^{m-1}\Psi(\theta_k,\xi_{k+1})\right\}
=\E\left[\exp\left\{\frac{1}{\sqrt m}\sum_{n=0}^{m-2}\Psi(\theta_k,\xi_{k+1})\right\}
\E_{m-1}\left[e^{\frac{1}{\sqrt m}\Psi(\theta_{m-1},\xi_m)}\right]\right].
\end{eqnarray*}
By Taylor expansion, we deduce that
\begin{align*}
\E_{m-1}\left[e^{\frac{1}{\sqrt m}\Psi(\theta_{m-1},\xi_m)}\right]&=\E_{m-1}\left[\sum_{n=0}^\infty\frac{1}{n!}\left(\frac{1}{\sqrt m}\Psi(\theta_{m-1},\xi_m)\right)^n\right]\\
&=1+\sum_{n=2}^\infty\E_{m-1}\left[\frac{1}{n!}\left(\frac{1}{\sqrt m}\Psi(\theta_{m-1},\xi_m)\right)^n\right]\\
&\le 1+\sum_{n=2}^\infty\E_{m-1}\left[\frac{1}{n!}\left(\frac{K}{\sqrt m}(1+|\xi_{m}|^2)\right)^n\right].
\end{align*}
For each element, we have
\begin{align*}
\E_{m-1}\left[\frac{1}{n!}\left(\frac{K}{\sqrt m}(1+|\xi_{m}|^2)\right)^n\right]
&\le \E_{m-1}\left[\frac{2^{n-1}}{n!}\left(\frac{K^n}{m^{\frac n2}}(1+|\xi_{m}|^{2n})\right)\right]\\
&\le \frac{2^{n-1}K^n}{n!m^{\frac n2}}(1+ d^{n-1}(2n-1)!!)\le \frac{(4Kd)^n}{m^{\frac n2}}.
\end{align*}
For small enough $\eta$ such that $\frac{4Kd }{\sqrt m}=4Kd \eta<1$, we have
\begin{align*}
\E_{m-1}\left[e^{\frac{1}{\sqrt m}\Psi(\theta_{m-1},\xi_m)}\right]
\le 1+\frac{\frac{(4Kd)^2}{m}}{1-\frac{4Kd}{\sqrt m}}=1+\frac{(4Kd)^2}{m-4Kd\sqrt m}.
\end{align*}
Inductively, we can get
\begin{align*}
\E \exp\left\{\frac{1}{\sqrt m}\sum_{n=0}^{m-1}\Psi(\theta_k,\xi_{k+1})\right\}
\le \left(1+\frac{(4Kd)^2}{m-4Kd\sqrt m}\right)^m\le C.
\end{align*}
\Ep

\begin{proof}[Proof of Lemma \ref{lem-R}]
We first consider the case for $\theta_0\sim\pi_\eta$. Recalling the definition of $\mathcal{R}_\eta$, we have
\begin{align*}
\PP(|\mathcal{R}_\eta|>x)\le\sum_{i=1}^6\PP(|\mathcal{R}_{\eta,i}|>\frac{x}{6}),
\end{align*}
and shall prove below that the following estimates hold:
\begin{eqnarray}
\PP(|\mathcal{R}_{\eta,1}|>x/6)&\le & Ce^{-c x\eta^{-\frac{1}{2}}}, \label{e:R-1} \\
\PP\left(\left\vert\mathcal{R}_{\eta,2}\right\vert>x/6\right) &\le & Ce^{-c x\eta^{-\frac 12}}, \label{e:R-2} \\
\PP\left(\left\vert\mathcal{R}_{\eta,3}\right\vert >x/6\right)& \le & Ce^{-cx\eta^{-\frac 32} }, \label{e:R-3}\\
\PP(|\mathcal{R}_{\eta,4}|>x/6) & \le &  C \left(e^{-c(x \eta^{-1})^{1/2}}1_{\{x<\eta^{-1}\}}+ e^{-c\eta^{-\frac 35}x^{\frac 25}}1_{\{x\ge \eta^{-1}\}}\right)\label{e:R-4},\\
\PP(|\mathcal{R}_{\eta,5}|>x/6)&\le & C e^{-cx^{\frac23}\eta^{-\frac 43} }, \label{e:R-5} \\
\PP(|\mathcal{R}_{\eta,6}|>x/6)&\le & C e^{-c\eta^{-2\bar\gamma}x^{\frac 23} }. \label{e:R-6}
\end{eqnarray}
where $\bar\gamma\in(0,1/4)$, $C$ and $c$ depends on $L, K_1, K_2, \sigma, |g(0)|^2$. Combining these estimates, we immediately conclude the proof.

 Let us show \eqref{e:R-1}-\eqref{e:R-6} below.
For \eqref{e:R-1}, by the Markov inequality and the fact that $\vp$ is uniformly bounded,
\begin{align*}
\PP(|\mathcal{R}_{\eta,1}|>\frac{x}{6})\le\E \exp\left\{|\vp(\theta_0)-\vp(\theta_m )|\right\} e^{-\frac x6\eta^{-\frac{1}{2}}} \le Ce^{-\frac x6\eta^{-\frac{1}{2}}}.
\end{align*}

For \eqref{e:R-2}, by the Markov inequality,
\begin{eqnarray*}
\PP\left(\mathcal{R}_{\eta,2}>\frac{x}{6}\right)
\le\E \exp\left\{\eta\sum_{k=0}^{m-1}\Ll\nabla^2 \vp(\theta_k ),(\sigma\xi_{k+1})(\sigma\xi_{k+1})^\mathrm{T}-\sigma\sigma^\mathrm{T}\Rr_\mathrm{HS}\right\} e^{- \frac x3\eta^{-\frac 12}}.
\end{eqnarray*}
Since
\begin{align*}
\E_k\left[\Ll\nabla^2 \vp(\theta_k ),(\sigma\xi_{k+1})(\sigma\xi_{k+1})^\mathrm{T}-\sigma\sigma^\mathrm{T}\Rr_\mathrm{HS}\right]=0
\end{align*}
and
\begin{eqnarray*}
\Ll\nabla^2 \vp(\theta_k ),(\sigma\xi_{k+1})(\sigma\xi_{k+1})^\mathrm{T}-\sigma\sigma^\mathrm{T}\Rr_\mathrm{HS}
&\le& \|\nabla^2 \vp\|\|(\sigma\xi_{k+1})(\sigma\xi_{k+1})^\mathrm{T}-\sigma\sigma^\mathrm{T}\|\\
&\le& K(1+|\xi_{k+1}|^2)
\end{eqnarray*}
with $K=\|\sigma\|^2  \|\nabla^2 \vp\|$,
by Lemma \ref{C5} with $\Psi(\theta_k,\xi_{k+1})=\Ll\nabla^2 \vp(\theta_k ),(\sigma\xi_{k+1})(\sigma\xi_{k+1})^\mathrm{T}-\sigma\sigma^\mathrm{T}\Rr_\mathrm{HS}$ therein, we have
\begin{eqnarray*}
\E \exp\left\{\eta\sum_{k=0}^{m-1}\Ll\nabla^2 \vp(\theta_k ),\sigma\sigma^\mathrm{T}-(\sigma\xi_{k+1})(\sigma\xi_{k+1})^\mathrm{T}\Rr_\mathrm{HS}\right\}
\le C.
\end{eqnarray*}
Therefore, $\PP\left(\mathcal{R}_{\eta,2}>\frac{x}{6}\right)
\le Ce^{-\frac x3 \eta^{-\frac 12}}$. Similarly, $\PP\left(\mathcal{R}_{\eta,2}<-\frac{x}{6}\right)
\le Ce^{-\frac x3\eta^{-\frac 12} }$, \eqref{e:R-2} is proved.

For \eqref{e:R-3} and $\gamma>0$, by the Markov inequality, we have
\begin{eqnarray*}
\PP  \left(\mathcal{R}_{\eta,3}>\frac{x}{6}\right)
\le \E \exp\left\{\sum_{k=0}^{m-1}\Ll (\eta \gamma)^{\frac 12}\left((\nabla^2\vp(\theta_k))^{\mathrm{T}}+\nabla^2\vp(\theta_k)\right)g(\theta_k),\sigma\xi_{k+1}\Rr\right\}e^{-\frac{x}{3}\sqrt\gamma\eta^{-\frac32}}.
\end{eqnarray*}
By Lemma \ref{C2} with $\Psi_1(\theta_k)=(\eta\gamma)^{\frac 12}\left((\nabla^2\vp(\theta_k))^{\mathrm{T}}+\nabla^2\vp(\theta_k)\right)g(\theta_k)$ and $\Psi_2=0$ therein, we have
\begin{align*}
&\E \exp\left\{\sum_{k=0}^{m-1}\Ll (\eta\gamma)^{\frac 12}\left((\nabla^2\vp(\theta_k))^{\mathrm{T}}+\nabla^2\vp(\theta_k)\right)g(\theta_k),\sigma\xi_{k+1}\Rr\right\}\\
\le& \left(\E \exp\left\{\sum_{k=0}^{m-1}2\eta\gamma |\left((\nabla^2\vp(\theta_k))^{\mathrm{T}}+\nabla^2\vp(\theta_k)\right)g(\theta_k)|^2\|\sigma\|^2\right\}\right)^{\frac 12} \\
\le& \left(\E \exp\left\{\sum_{k=0}^{m-1}C \gamma\eta |g(\theta_k)|^2\right\} \right)^{\frac 12}.
\end{align*}
Choosing $\gamma$ small enough such that $C \gamma\le\gamma_0$ in Lemma \ref{lem-g}, and combining the previous two relations with (\ref{R_1,22}), we obtain
$\PP\left(\mathcal{R}_{\eta,3} > \frac{x}{6} \right)
\le C\exp\{c(\eta^{-1}-x\eta^{-\frac 32})\}$. By the same argument, we
obtain the same bound for $\PP\left(\mathcal{R}_{\eta,3}<-\frac{x}{6} \right)$. Hence for $x>\eta^{\frac12}$, we have
\begin{align*}
\PP\left(\left\vert\mathcal{R}_{\eta,3}\right\vert > \frac{x}{6} \right)
\le C\exp\left\{-c\eta^{-1}(x\eta^{-\frac 12}-1)\right\}
\le C\exp\left\{-cx\eta^{-\frac 32}\right\}.
\end{align*}

For \eqref{e:R-4}, we have
\begin{eqnarray*}
&&\PP(\mathcal{R}_{\eta,4}>\frac{x}{6})\\
&=&\PP\left(\sum_{k=0}^{m-1}\int_0^1\sum_{i_1,i_2,i_3=1}^d\nabla^3_{i_1,i_2,i_3}\vp(\theta_k+t\Delta\theta_k) (\sigma\xi_{k+1})_{i_1}(\sigma\xi_{k+1})_{i_2}(\sigma\xi_{k+1})_{i_3} \dif t>\eta^{-2}x\right)\\
&\le&\PP\left(\sum_{k=0}^{m-1}\int_0^1\sum_{i_1,i_2,\atop i_3=1}^d
\left(\nabla^3_{i_1,i_2,i_3}\vp(\theta_k+t\Delta\theta_k)-\nabla^3_{i_1,i_2,i_3}\vp(\theta_k) \right) (\sigma\xi_{k+1})_{i_1}(\sigma\xi_{k+1})_{i_2}(\sigma\xi_{k+1})_{i_3} \dif t>\frac{\eta^{-2}x}{2}\right)\\
&&+\PP\left(\sum_{k=0}^{m-1} \sum_{i_1,i_2,i_3=1}^d\nabla^3_{i_1,i_2,i_3}\vp(\theta_k) (\sigma\xi_{k+1})_{i_1}(\sigma\xi_{k+1})_{i_2}(\sigma\xi_{k+1})_{i_3}  >\frac{\eta^{-2}x}{2}\right)\\
&:=&\mathcal{R}_{\eta,4,1}+\mathcal{R}_{\eta,4,2}.
\end{eqnarray*}
For $\mathcal{R}_{\eta,4,1}$, applying Taylor expansion to the function $\nabla^3_{i_1,i_2,i_3}\vp$, we get
\begin{eqnarray*}
&&\mathcal{R}_{\eta,4,1}\\
&=&\PP\left(\sum_{k=0}^{m-1}\int_0^1\int_0^1\sum_{i_1,i_2,\atop i_3,i_4=1}^d
\nabla^4_{i_1,i_2, i_3,i_4}\vp(\theta_k+tt'\Delta\theta_k)(t\Delta\theta_k)_{i_4} (\sigma\xi_{k+1})_{i_1}(\sigma\xi_{k+1})_{i_2}(\sigma\xi_{k+1})_{i_3} \dif t' \dif t>\frac{\eta^{-2}x}{2}\right)\\
&\le&\PP\left(\sum_{k=0}^{m-1}
|\eta g(\theta_k)+\sqrt \eta\sigma\xi_{k+1}| |\sigma\xi_{k+1}|^3> c\eta^{-2}x\right),
\end{eqnarray*}
where the last inequality follows the boundedness of $\|\nabla^4\vp\|$ and Cauchy's inequality. Then we have
\begin{equation}\label{e-R_4,1}
\mathcal{R}_{\eta,4,1}\le\PP\left(\sum_{k=0}^{m-1}
|\eta g(\theta_k)| |\sigma\xi_{k+1}|^3> c\eta^{-2}x\right)
+\PP\left(\sum_{k=0}^{m-1}
\sqrt \eta|\sigma\xi_{k+1}|^4> c\eta^{-2}x\right).
\end{equation}
For the first probability, denoting $A=\{|\xi_i| \le c \eta^{-\delta} x^{\delta}, i=1,...,m\}$ with $\delta \in (0,1/3)$, we have
\begin{eqnarray*}
&&\PP\left(\sum_{k=0}^{m-1}
|\eta g(\theta_k)| |\sigma\xi_{k+1}|^3>c\eta^{-2}x\right)\\
&\le& \PP\left(\sum_{k=0}^{m-1}
|\eta g(\theta_k)| |\sigma\xi_{k+1}|^3>c\eta^{-2}x, A\right)
+\PP\left(A^c\right)\\
&\le& \PP\left(\sum_{k=0}^{m-1}
\eta|g(\theta_k)|^2>C\eta^{-1} (x \eta^{-1})^{2-6\delta}\right)
+\sum_{k=0}^{m-1}\PP\left(|\xi_{k+1}|>c(\eta^{-1}x)^{\delta}\right),
\end{eqnarray*}
where the last inequality is because the following relation holds in the set $A$:
\Bes
\sum_{k=0}^{m-1}
|\eta g(\theta_k)| |\sigma\xi_{k+1}|^3 \le C\eta^{1-3 \delta}x^{3 \delta} \sum_{k=0}^{m-1}
|g(\theta_k)| \le C\eta^{-3 \delta}x^{3 \delta} \left(\sum_{k=0}^{m-1}
|g(\theta_k)|^2\right)^{\frac12}.
\Ees
For the first term, \eqref{e:R-1-22} yields
\begin{eqnarray*}
\PP\left(\sum_{k=0}^{m-1}
\eta|g(\theta_k)|^2>C\eta^{-1} (x \eta^{-1})^{2-6\delta}\right)
\le e^{-C\eta^{-1}((x\eta^{-1})^{2-6 \delta}-1)}.
\end{eqnarray*}
For the second term, the tail probability estimate of gaussian distribution implies
\begin{eqnarray*}
\sum_{k=0}^{m-1}\PP\left(|\xi_{k+1}|>c(\eta^{-1}x)^{\delta}\right)\le Cme^{-c(x \eta^{-1})^{2 \delta}}.
\end{eqnarray*}
Hence we have
\begin{eqnarray}
\PP\left(\sum_{k=0}^{m-1}
|\eta g(\theta_k)| |\sigma\xi_{k+1}|^3>c\eta^{-2}x\right)
&\le& e^{-C\eta^{-1}((x\eta^{-1})^{2-6 \delta}-1)}+C m e^{-c(x \eta^{-1})^{2 \delta}}\nonumber\\
&\le& C \exp\left(-c(x \eta^{-1})^{1/2}\right)\label{e-R_4,1,1}
\end{eqnarray}
as $x\ge\eta^{\frac12}$ and $\delta=1/4$.

For the second probability of (\ref{e-R_4,1}), denoting $\hat\xi_{k+1}=\xi_{k+1}1_{\{|\xi_{k+1}|<\eta^{-\delta}x^{\delta'}\}}$ and $ \check{\xi}_{k+1}=\xi_{k+1}1_{\{|\xi_{k+1}|\ge\eta^{-\delta}x^{\delta'}\}}$, we have
\begin{eqnarray*}
\PP\left(\sum_{k=0}^{m-1}\sqrt \eta|\sigma\xi_{k+1}|^4> c\eta^{-2}x\right)
&\le& \PP\left(\sum_{k=0}^{m-1}|\xi_{k+1}|^4> C\eta^{-\frac 52}x\right)\\
&\le& \PP\left(\sum_{k=0}^{m-1}|\hat{\xi}_{k+1}|^4> C\eta^{-\frac 52}x\right)+\PP\left(\sum_{k=0}^{m-1}|\check{\xi}_{k+1}|^4> C\eta^{-\frac 52}x\right).
\end{eqnarray*}
Let us bound the two terms on the right hand side. For the first one,
\begin{eqnarray*}
\PP\left(\sum_{k=0}^{m-1}|\hat{\xi}_{k+1}|^4> C\eta^{-\frac 52}x\right)
&=&\PP\left(\sum_{k=0}^{m-1}\left(|\hat{\xi}_{k+1}|^4-\E[|\hat{\xi}_{k+1}|^4]\right)> C\eta^{-\frac 52}\left(x-\eta^{\frac12}\E[|\hat{\xi}_{1}|^4]\right)\right) \\
&\le& \exp\left\{-\frac{C\eta^{-5}{\left(x-\eta^{\frac12}\E[|\hat{\xi}_{k+1}|^4]\right)^2}}{\eta^{-2}(c\eta^{-4\delta}x^{4\delta'})^2}\right\}\\
&\le& \exp\left\{-C\eta^{-3+8\delta}x^{2-8\delta'}\right\}
\end{eqnarray*}
with $x\ge c\eta^{\frac 12}$, where the first inequality is by Hoeffding's inequality and the bound $||\hat{\xi}_{k+1}|^4-\E[|\hat{\xi}_{k+1}|^4]|\le c \eta^{-4\delta}x^{4\delta'}$.
For the second term, it follows from Chebyshev's inequality that
\begin{eqnarray}
\PP\left(\sum_{k=0}^{m-1}|\check{\xi}_{k+1}|^4> C\eta^{-\frac 52}x\right)
&\le&\frac{C}{x}\eta^{\frac 52}\E\left[\sum_{k=0}^{m-1}|\check{\xi}_{k+1}|^4\right] \nonumber\\
&=&\frac{C}{x}\eta^{\frac 52}m\int_{|y| \ge \eta^{-\delta}x^{\delta'}} \frac{1}{(2\pi)^{d/2}}|y|^4 e^{-\frac 12 |y|^2}dy\nonumber\\
&\le& C x^{4\delta'-1} \eta^{\frac 12-4\delta}e^{-c\eta^{-2\delta}x^{2\delta'}}. \label{e:R41-1}
\end{eqnarray}
Combining the estimates for the two probabilities of (\ref{e-R_4,1}) and taking $\delta=\frac {3}{10}$, $\delta'=\frac 15$, we obtain
\begin{eqnarray*}
\mathcal{R}_{\eta,4,1}&\le& C \left(\exp\left\{-c(x \eta^{-1})^{1/2}\right\}+\exp\left\{-c\eta^{-\frac 35}x^{\frac 25}\right\}\right)\\
&\le& C \left(\exp\left\{-c(x \eta^{-1})^{1/2}\right\}1_{\{x<\eta^{-1}\}}+ \exp\left\{-c\eta^{-\frac 35}x^{\frac 25}\right\}1_{\{x\ge \eta^{-1}\}}\right)
\end{eqnarray*}
as $x\ge c\eta^{\frac 12}$.

  We now estimate $\mathcal{R}_{\eta,4,2}$. For $\delta>0$,  denoting $\hat\xi_{k+1}=\xi_{k+1}1_{\{|\xi_{k+1}|< c\eta^{-\delta}x^{\delta'}\}}$ and $ \check{\xi}_{k+1}=\xi_{k+1}1_{\{|\xi_{k+1}|\ge c\eta^{-\delta}x^{\delta'}\}}$, we have
\begin{align}
\mathcal{R}_{\eta,4,2}\le &\PP\left(\sum_{k=0}^{m-1} \sum_{i_1,i_2,i_3=1}^d\nabla^3_{i_1,i_2,i_3}\vp(\theta_k) (\sigma\hat{\xi}_{k+1})_{i_1}(\sigma\hat{\xi}_{k+1})_{i_2}(\sigma\hat{\xi}_{k+1})_{i_3}  >\frac{\eta^{-2}x}{2}\right) \nonumber \\
&+\PP\left(\sum_{k=0}^{m-1} \sum_{i_1,i_2,i_3=1}^d\nabla^3_{i_1,i_2,i_3}\vp(\theta_k) (\sigma\check{\xi}_{k+1})_{i_1}(\sigma\check{\xi}_{k+1})_{i_2}(\sigma\check{\xi}_{k+1})_{i_3}>\frac{\eta^{-2}x}{2}\right). \label{e:R42}
\end{align}
Let us bound the above two probabilities. For the first one, let $\lambda>0$, we have
\begin{align*}
&\PP\left(\sum_{k=0}^{m-1} \sum_{i_1,i_2,i_3=1}^d\nabla^3_{i_1,i_2,i_3}\vp(\theta_k) (\sigma\hat{\xi}_{k+1})_{i_1}(\sigma\hat{\xi}_{k+1})_{i_2}(\sigma\hat{\xi}_{k+1})_{i_3}>\frac{\eta^{-2}x}{2}\right)\\
\le&\E \exp\left\{\sum_{k=0}^{m-1} \lambda\sum_{i_1,i_2,i_3=1}^d\nabla^3_{i_1,i_2,i_3}\vp(\theta_k) (\sigma\hat{\xi}_{k+1})_{i_1}(\sigma\hat{\xi}_{k+1})_{i_2}(\sigma\hat{\xi}_{k+1})_{i_3}\right\} e^{- \frac x2\lambda\eta^{-2}}.
\end{align*}
Since
\begin{align*}
\E_k\left[\sum_{i_1,i_2,i_3=1}^d\nabla^3_{i_1,i_2,i_3}\vp(\theta_k) (\sigma\hat{\xi}_{k+1})_{i_1}(\sigma\hat{\xi}_{k+1})_{i_2}(\sigma\hat{\xi}_{k+1})_{i_3}\right]=0,
\end{align*}
and
\begin{align*}
\left|\sum_{i_1,i_2,i_3=1}^d\nabla^3_{i_1,i_2,i_3}\vp(\theta_k) (\sigma\hat{\xi}_{k+1})_{i_1}(\sigma\hat{\xi}_{k+1})_{i_2}(\sigma\hat{\xi}_{k+1})_{i_3}\right|
\le\|\nabla^3\vp\||\sigma\hat{\xi}_{k+1}|^3 \le C\eta^{-3\delta}x^{3\delta'},
\end{align*}
Hoeffding's lemma (cf.\citet[Lemma 2.6]{Massart2003}) gives
\begin{align*}
\E_k \exp\left\{\lambda\sum_{i_1,i_2,i_3=1}^d\nabla^3_{i_1,i_2,i_3}\vp(\theta_k) (\sigma\hat{\xi}_{k+1})_{i_1}(\sigma\hat{\xi}_{k+1})_{i_2}(\sigma\hat{\xi}_{k+1})_{i_3} \right\}\le e^{C\lambda^2\eta^{-6\delta}x^{6\delta'}}.
\end{align*}
Using conditional expectation inductively, we obtain
\begin{align*}
& \E \exp\left\{\sum_{k=0}^{m-1} \lambda\sum_{i_1,i_2,i_3=1}^d\nabla^3_{i_1,i_2,i_3}\vp(\theta_k) (\sigma\hat{\xi}_{k+1})_{i_1}(\sigma\hat{\xi}_{k+1})_{i_2}(\sigma\hat{\xi}_{k+1})_{i_3}\right\}  \\
 =&\E\E_{m-1} \exp\left\{\sum_{k=0}^{m-1} \lambda\sum_{i_1,i_2,i_3=1}^d\nabla^3_{i_1,i_2,i_3}\vp(\theta_k) (\sigma\hat{\xi}_{k+1})_{i_1}(\sigma\hat{\xi}_{k+1})_{i_2}(\sigma\hat{\xi}_{k+1})_{i_3}\right\}  \\
 \le& e^{C\lambda^2\eta^{-6\delta}x^{6\delta'}} \E \exp\left\{\sum_{k=0}^{m-2} \lambda\sum_{i_1,i_2,i_3=1}^d\nabla^3_{i_1,i_2,i_3}\vp(\theta_k) (\sigma\hat{\xi}_{k+1})_{i_1}(\sigma\hat{\xi}_{k+1})_{i_2}(\sigma\hat{\xi}_{k+1})_{i_3}\right\}  \\
 \le& ... \le e^{C\lambda^2\eta^{-6\delta-2}x^{6\delta'}}.
\end{align*}
Taking $\lambda=\frac {1}{4C}\eta^{6\delta}x^{1-6\delta'}$, hence we have
\begin{eqnarray*}
 \PP\left(\sum_{k=0}^{m-1} \sum_{i_1,i_2,i_3=1}^d\nabla^3_{i_1,i_2,i_3}\vp(\theta_k) (\sigma\hat{\xi}_{k+1})_{i_1}(\sigma\hat{\xi}_{k+1})_{i_2}(\sigma\hat{\xi}_{k+1})_{i_3} >\frac{\eta^{-2}x}{2}\right)
\le e^{-Cx^{2-6\delta'}\eta^{6\delta-2}}.
\end{eqnarray*}
For the second probability of \eqref{e:R42}, we have
\begin{eqnarray*}
\PP\left(\sum_{k=0}^{m-1} \sum_{i_1,i_2,i_3=1}^d\nabla^3_{i_1,i_2,i_3}\vp(\theta_k) (\sigma\check{\xi}_{k+1})_{i_1}(\sigma\check{\xi}_{k+1})_{i_2}(\sigma\check{\xi}_{k+1})_{i_3} >\frac{\eta^{-2}x}{2}\right)
\le\PP\left(\sum_{k=0}^{m-1} |\check{\xi}_{k+1}|^3  >C\eta^{-2}x\right),
\end{eqnarray*}
which, together with the Markov inequality and a similar argument as in \eqref{e:R41-1}, yields
\begin{align*}
\PP\left(\sum_{k=0}^{m-1}|\check{\xi}_{k+1}|^3>C\eta^{-2}x\right)
\le\frac{C\eta^2}{x}\sum_{k=0}^{m-1}\E\left[|\check{\xi}_{k+1}|^3\right]
\le &Ce^{-c \eta^{-2\delta}x^{2\delta'}}
\end{align*}
as $x^{\delta'}\ge\eta^{\frac\delta2}$. Choosing $\delta=\delta'=\frac {1}{4}$, one gets
$$
\mathcal{R}_{\eta,4,2}\le C e^{-cx^{\frac 12}\eta^{-\frac 12}}.
$$
Combining the estimates of $\mathcal{R}_{\eta,4,1}$ and $\mathcal{R}_{\eta,4,2}$, we obtain
\begin{eqnarray*}
\PP(\mathcal{R}_{\eta,4}\ge\frac{x}{6}) &\le& \mathcal{R}_{\eta,4,1}+\mathcal{R}_{\eta,4,2}\\
&\le& C \left(\exp\left\{-c(x \eta^{-1})^{1/2}\right\}1_{\{x<\eta^{-1}\}}+ \exp\left\{-c\eta^{-\frac 35}x^{\frac 25}\right\}1_{\{x\ge \eta^{-1}\}}\right)
\end{eqnarray*}
as $x\ge c\eta^{\frac 12}$.

For (\ref{e:R-5}), the boundedness of $\|\nabla^2\vp\|$ and $\|\nabla^3\vp\|$  implies
\begin{eqnarray*}
\PP(|\mathcal{R}_{\eta,5}|>\frac x6)&\le&
\PP\left(\eta^{\frac{5}{2}}\sum_{k=0}^{m-1}|g(\theta_k )|^2>Cx\right)+\PP\left(\eta^{\frac 72}\sum_{k=0}^{m-1}|g(\theta_k )|^3>Cx\right)\\
&\le&\PP\left(\eta\sum_{k=0}^{m-1}|g(\theta_k )|^2>Cx\eta^{-\frac 32}\right)+\PP\left(\sum_{k=0}^{m-1}\eta|g(\theta_k )|^2>Cx^{\frac23}\eta^{-\frac43}\right)\\
&\le&C\left(\exp\left\{-c\eta^{-1}(x\eta^{-\frac 12}-1)\right\}+\exp\left\{-c\eta^{-1}(x^{\frac23}\eta^{-\frac 13}-1)\right\}\right)\\
&\le& C\exp\left\{-cx^{\frac23}\eta^{-\frac 43}\right\}
\end{eqnarray*}
as $x\ge c\eta^{1/2}$, where the third line is by \eqref{e:R-1-22}.

For (\ref{e:R-6}), by the boundedness of $\|\nabla^3\vp\|$, one has
\begin{eqnarray*}
\PP(|\mathcal{R}_{\eta,6}|\ge \frac{x}{6})
&\le&\PP\left(\eta^{\frac 52}\sum_{k=0}^{m-1}\left(|g(\theta_k)||\sigma\xi_{k+1}|^2+\sqrt\eta|g(\theta_k)|^2|\sigma\xi_{k+1}|\right)\ge Cx\right)\\
&\le&\PP\left(\sum_{k=0}^{m-1}|g(\theta_k)||\sigma\xi_{k+1}|^2\ge Cx\eta^{-\frac 52}\right)
+\PP\left(\sum_{k=0}^{m-1}|g(\theta_k)|^2|\sigma\xi_{k+1}|\ge Cx\eta^{-3}\right).
\end{eqnarray*}
For the first probability, denoting $A=\{|\xi_i|\le c\eta^{-\bar\gamma}x^{\bar\gamma'},i=1,...,m\}$, we have
\begin{eqnarray*}
&&\PP\left(\sum_{k=0}^{m-1}|g(\theta_k)||\sigma\xi_{k+1}|^2\ge Cx\eta^{-\frac 52}\right)\\
&\le& \PP\left(\sum_{k=0}^{m-1}|g(\theta_k)||\sigma\xi_{k+1}|^2\ge Cx\eta^{-\frac 52},A\right)+\PP(A^c)\\
&\le& \PP\left(\sum_{k=0}^{m-1}\eta|g(\theta_k)|^2\ge Cx^{2-4\bar\gamma'}\eta^{-2+4\bar\gamma}\right) +\sum_{k=0}^{m-1}\PP(|\xi_{k+1}|>c\eta^{-\bar\gamma}x^{\bar\gamma'}).
\end{eqnarray*}
Similar with the calculation of (\ref{e-R_4,1,1}), we have
\begin{eqnarray*}
\PP\left(\sum_{k=0}^{m-1}|g(\theta_k)||\sigma\xi_{k+1}|^2\ge Cx\eta^{-\frac 52}\right)
&\le& C\left(e^{-c\eta^{-1}(x^{2-4\bar\gamma'}\eta^{-1+4\bar\gamma}-1)}+e^{-c\eta^{-2\bar\gamma}x^{2\bar\gamma'}}\right)\\
&\le& C\exp\left\{-c\eta^{-2\bar\gamma}x^{\frac 23}\right\}
\end{eqnarray*}
as $\bar\gamma'=1/3$, $\bar\gamma<\frac 14$ and $x^{\frac 23}\ge\max\{c \eta^{1-4\bar\gamma},\eta^{\bar\gamma}\}$,  that is, $x\ge c\max\{\eta^{\frac32-6\bar\gamma},\eta^{\frac32\bar\gamma}\}$.
For the second probability, we use the same method with $A=\{|\xi_i|\le c\eta^{-\delta}x^{\delta'},i=1,...,m\}$  and get
\begin{eqnarray*}
&&\PP\left(\sum_{k=0}^{m-1}|g(\theta_k)|^2|\sigma\xi_{k+1}|\ge Cx\eta^{-3}\right)\\
&\le& \PP\left(\sum_{k=0}^{m-1}\eta|g(\theta_k)|^2\ge Cx^{1-\delta'}\eta^{-2+\delta}\right)+\sum_{k=0}^{m-1}\PP\left(|\sigma\xi_{k+1}|\ge x^{\delta'}\eta^{-\delta}\right)\\
&\le& C\left(\exp\{-c\eta^{-1}(x^{1-\delta'}\eta^{-1+\delta}-1)\}+m\exp\{-cx^{2\delta'}\eta^{-2\delta}\}\right)\\
&\le& C\exp\{-cx^{\frac23}\eta^{-\frac43}\}
\end{eqnarray*}
as $\delta'=1/3$, $\delta=2/3$ and $x\ge c\eta^{1/2}$.
Hence, we have
\begin{align*}
\PP(|\mathcal{R}_{\eta,6}|\ge \frac{x}{6})\le C\exp\{-c\eta^{-2\bar\gamma}x^{\frac 23}\}
\end{align*}
with $x\ge c\max\{\eta^{\frac 32-6\bar\gamma},\eta^{\frac32\bar\gamma},\eta^{\frac 12}\}$ as $\bar\gamma<\frac 14$.
Combining (\ref{e:R-1}-\ref{e:R-6}), we have
\begin{eqnarray*}
\PP(|\mathcal{R}_\eta|>x)\le& C \left(e^{-c(x \eta^{-1})^{1/2}}1_{\{x<\eta^{-1}\}}+ e^{-c\eta^{-\frac 35}x^{\frac 25}}1_{\{x\ge \eta^{-1}\}}+e^{-c\eta^{-2\bar\gamma}x^{\frac23}}\right),
\end{eqnarray*}
as  $x\ge c\max\{\eta^{\frac 32-6\bar\gamma},\eta^{\frac32\bar\gamma},\eta^{\frac 12}\}$ where $0<\bar\gamma<\frac{1}{4}$.

The case without condition $\theta_0\sim\pi_\eta$ can be estimate by a similar way following \eqref{R_1,22-0} and \eqref{e:R-1-22-0}.

\end{proof}

\end{appendix}

\section*{Acknowledgements}
We would like to gratefully thank Professors Fuqing Gao and Feng-Yu Wang for very helpful discussions.  We also thank two anonymous referees and the AE for their valuable comments which have improved the manuscript considerably. LX is supported in part by Macao S.A.R grant FDCT  0090/2019/A2 and University of Macau grant  MYRG2018-00133-FST.

\bibliographystyle{imsart-nameyear} 
\bibliographystyle{imsart-nameyear} 

\begin{thebibliography}{40}

\bibitem[\protect\citeauthoryear{Applebaum}{2009}]{applebaum2009levy}
\begin{bbook}[author]
\bauthor{\bsnm{Applebaum},~\bfnm{David}\binits{D.}}
(\byear{2009}).
\btitle{L{\'e}vy processes and stochastic calculus}.
\bpublisher{Cambridge university press}.
\end{bbook}
\endbibitem

\bibitem[\protect\citeauthoryear{Bao, Huang and Yuan}{2019}]{Bao19}
\begin{barticle}[author]
\bauthor{\bsnm{Bao},~\bfnm{Jianhai}\binits{J.}},
  \bauthor{\bsnm{Huang},~\bfnm{Xing}\binits{X.}} \AND
  \bauthor{\bsnm{Yuan},~\bfnm{Chenggui}\binits{C.}}
(\byear{2019}).
\btitle{New regularity of kolmogorov equation and application on approximation
  of semi-linear spdes with H{\"o}lder continuous drifts}.
\bjournal{Communications on Pure \& Applied Analysis}
\bvolume{18}
\bpages{341--361}.
\end{barticle}
\endbibitem

\bibitem[\protect\citeauthoryear{Bao and Shao}{2018}]{Bao18}
\begin{barticle}[author]
\bauthor{\bsnm{Bao},~\bfnm{Jianhai}\binits{J.}} \AND
  \bauthor{\bsnm{Shao},~\bfnm{Jinghai}\binits{J.}}
(\byear{2018}).
\btitle{Weak convergence of path-dependent SDEs with irregular coefficients}.
\bjournal{arXiv preprint arXiv:1809.03088}.
\end{barticle}
\endbibitem

\bibitem[\protect\citeauthoryear{Bao and Yuan}{2013}]{Bao12}
\begin{barticle}[author]
\bauthor{\bsnm{Bao},~\bfnm{Jianhai}\binits{J.}} \AND
  \bauthor{\bsnm{Yuan},~\bfnm{Chenggui}\binits{C.}}
(\byear{2013}).
\btitle{Convergence rate of EM scheme for SDDEs}.
\bjournal{Proceedings of the American Mathematical Society}
\bvolume{141}
\bpages{3231--3243}.
\end{barticle}
\endbibitem

\bibitem[\protect\citeauthoryear{Barlow, Jacka and
  Yor}{1986}]{barlow1986inequalities}
\begin{barticle}[author]
\bauthor{\bsnm{Barlow},~\bfnm{MT}\binits{M.}},
  \bauthor{\bsnm{Jacka},~\bfnm{SD}\binits{S.}} \AND
  \bauthor{\bsnm{Yor},~\bfnm{M}\binits{M.}}
(\byear{1986}).
\btitle{Inequalities for a pair of processes stopped at a random time}.
\bjournal{Proceedings of the London Mathematical Society}
\bvolume{3}
\bpages{142--172}.
\end{barticle}
\endbibitem

\bibitem[\protect\citeauthoryear{Chatterji et~al.}{2020}]{JY3}
\begin{binproceedings}[author]
\bauthor{\bsnm{Chatterji},~\bfnm{Niladri}\binits{N.}},
  \bauthor{\bsnm{Diakonikolas},~\bfnm{Jelena}\binits{J.}},
  \bauthor{\bsnm{Jordan},~\bfnm{Michael~I.}\binits{M.~I.}} \AND
  \bauthor{\bsnm{Bartlett},~\bfnm{Peter}\binits{P.}}
(\byear{2020}).
\btitle{Langevin monte carlo without smoothness}.
In \bbooktitle{International Conference on Artificial Intelligence and
  Statistics}
\bpages{1716--1726}.
\bpublisher{PMLR}.
\end{binproceedings}
\endbibitem

\bibitem[\protect\citeauthoryear{Chen, Fang and Shao}{2013}]{CFS13}
\begin{barticle}[author]
\bauthor{\bsnm{Chen},~\bfnm{Louis H.~Y.}\binits{L.~H.~Y.}},
  \bauthor{\bsnm{Fang},~\bfnm{Xiao}\binits{X.}} \AND
  \bauthor{\bsnm{Shao},~\bfnm{Qi-Man}\binits{Q.-M.}}
(\byear{2013}).
\btitle{From Stein identities to moderate deviations}.
\bjournal{The Annals of Probability}
\bvolume{41}
\bpages{262--293}.
\end{barticle}
\endbibitem

\bibitem[\protect\citeauthoryear{Chen et~al.}{2016}]{CSWX16}
\begin{barticle}[author]
\bauthor{\bsnm{Chen},~\bfnm{Xiaohong}\binits{X.}},
  \bauthor{\bsnm{Shao},~\bfnm{Qi-Man}\binits{Q.-M.}},
  \bauthor{\bsnm{Wu},~\bfnm{Wei~Biao}\binits{W.~B.}} \AND
  \bauthor{\bsnm{Xu},~\bfnm{Lihu}\binits{L.}}
(\byear{2016}).
\btitle{Self-normalized Cram{\'e}r-type moderate deviations under dependence}.
\bjournal{The Annals of Statistics}
\bvolume{44}
\bpages{1593--1617}.
\end{barticle}
\endbibitem

\bibitem[\protect\citeauthoryear{Dalalyan}{2017}]{Dal17}
\begin{barticle}[author]
\bauthor{\bsnm{Dalalyan},~\bfnm{Arnak~S.}\binits{A.~S.}}
(\byear{2017}).
\btitle{Theoretical guarantees for approximate sampling from smooth and
  log-concave densities}.
\bjournal{Journal of the Royal Statistical Society Series B}
\bvolume{79}
\bpages{651--676}.
\end{barticle}
\endbibitem

\bibitem[\protect\citeauthoryear{Dedecker and
  Gou{\"e}zel}{2015}]{dedecker2015subgaussian}
\begin{barticle}[author]
\bauthor{\bsnm{Dedecker},~\bfnm{J{\'e}r{\^o}me}\binits{J.}} \AND
  \bauthor{\bsnm{Gou{\"e}zel},~\bfnm{S{\'e}bastien}\binits{S.}}
(\byear{2015}).
\btitle{Subgaussian concentration inequalities for geometrically ergodic Markov
  chains}.
\bjournal{Electronic Communications in Probability}
\bvolume{20}.
\end{barticle}
\endbibitem

\bibitem[\protect\citeauthoryear{Del~Moral, Hu and Wu}{2015}]{Moral15}
\begin{barticle}[author]
\bauthor{\bsnm{Del~Moral},~\bfnm{Pierre}\binits{P.}},
  \bauthor{\bsnm{Hu},~\bfnm{Shulan}\binits{S.}} \AND
  \bauthor{\bsnm{Wu},~\bfnm{Liming}\binits{L.}}
(\byear{2015}).
\btitle{Moderate deviations for interacting processes}.
\bjournal{Statistica Sinica}
\bvolume{25}
\bpages{921--951}.
\end{barticle}
\endbibitem

\bibitem[\protect\citeauthoryear{Dupuis and Johnson}{2017}]{DuJo17}
\begin{barticle}[author]
\bauthor{\bsnm{Dupuis},~\bfnm{Paul}\binits{P.}} \AND
  \bauthor{\bsnm{Johnson},~\bfnm{Dane}\binits{D.}}
(\byear{2017}).
\btitle{Moderate deviations-based importance sampling for stochastic recursive
  equations}.
\bjournal{Advances in Applied Probability}
\bvolume{49}
\bpages{981--1010}.
\end{barticle}
\endbibitem

\bibitem[\protect\citeauthoryear{Durmus and Moulines}{2017}]{JY1}
\begin{barticle}[author]
\bauthor{\bsnm{Durmus},~\bfnm{Alain}\binits{A.}} \AND
  \bauthor{\bsnm{Moulines},~\bfnm{{\'E}ric}\binits{{\'E}.}}
(\byear{2017}).
\btitle{Nonasymptotic convergence analysis for the unadjusted Langevin
  algorithm}.
\bjournal{The Annals of Applied Probability}
\bvolume{27}
\bpages{1551--1587}.
\end{barticle}
\endbibitem

\bibitem[\protect\citeauthoryear{Durmus and Moulines}{2019}]{JY2}
\begin{barticle}[author]
\bauthor{\bsnm{Durmus},~\bfnm{Alain}\binits{A.}} \AND
  \bauthor{\bsnm{Moulines},~\bfnm{{\'E}ric}\binits{{\'E}.}}
(\byear{2019}).
\btitle{High-dimensional Bayesian inference via the unadjusted Langevin
  algorithm}.
\bjournal{Bernoulli}
\bvolume{25}
\bpages{2854--2882}.
\end{barticle}
\endbibitem

\bibitem[\protect\citeauthoryear{Fan}{2020}]{fan2020cramer}
\begin{barticle}[author]
\bauthor{\bsnm{Fan},~\bfnm{Xiequan}\binits{X.}}
(\byear{2020}).
\btitle{Cram{\'e}r type moderate deviations for self-normalized $\psi$-mixing
  sequences}.
\bjournal{Journal of Mathematical Analysis and Applications}
\bvolume{486}
\bpages{123902}.
\end{barticle}
\endbibitem

\bibitem[\protect\citeauthoryear{Fan et~al.}{2019}]{Fan2019Self}
\begin{barticle}[author]
\bauthor{\bsnm{Fan},~\bfnm{Xiequan}\binits{X.}},
  \bauthor{\bsnm{Grama},~\bfnm{Ion}\binits{I.}},
  \bauthor{\bsnm{Liu},~\bfnm{Quansheng}\binits{Q.}} \AND
  \bauthor{\bsnm{Shao},~\bfnm{Qi-Man}\binits{Q.-M.}}
(\byear{2019}).
\btitle{Self-normalized Cram{\'e}r type moderate deviations for martingales}.
\bjournal{Bernoulli}
\bvolume{25}
\bpages{2793--2823}.
\end{barticle}
\endbibitem

\bibitem[\protect\citeauthoryear{Fan et~al.}{2020}]{Fan2020}
\begin{barticle}[author]
\bauthor{\bsnm{Fan},~\bfnm{Xiequan}\binits{X.}},
  \bauthor{\bsnm{Grama},~\bfnm{Ion}\binits{I.}},
  \bauthor{\bsnm{Liu},~\bfnm{Quansheng}\binits{Q.}} \AND
  \bauthor{\bsnm{Shao},~\bfnm{Qi-Man}\binits{Q.-M.}}
(\byear{2020}).
\btitle{Self-normalized Cram{\'e}r type moderate deviations for stationary
  sequences and applications}.
\bjournal{Stochastic Processes and their Applications}
\bvolume{130}
\bpages{5124--5148}.
\end{barticle}
\endbibitem

\bibitem[\protect\citeauthoryear{Fang, Luo and Shao}{2020}]{FLS20}
\begin{barticle}[author]
\bauthor{\bsnm{Fang},~\bfnm{Xiao}\binits{X.}},
  \bauthor{\bsnm{Luo},~\bfnm{Li}\binits{L.}} \AND
  \bauthor{\bsnm{Shao},~\bfnm{Qi-Man}\binits{Q.-M.}}
(\byear{2020}).
\btitle{A refined Cram{\'e}r-type moderate deviation for sums of local
  statistics}.
\bjournal{Bernoulli}
\bvolume{26}
\bpages{2319--2352}.
\end{barticle}
\endbibitem

\bibitem[\protect\citeauthoryear{Fang, Shao and Xu}{2019}]{FSX19}
\begin{barticle}[author]
\bauthor{\bsnm{Fang},~\bfnm{Xiao}\binits{X.}},
  \bauthor{\bsnm{Shao},~\bfnm{Qi-Man}\binits{Q.-M.}} \AND
  \bauthor{\bsnm{Xu},~\bfnm{Lihu}\binits{L.}}
(\byear{2019}).
\btitle{Multivariate approximations in Wasserstein distance by Stein's method
  and Bismut's formula}.
\bjournal{Probability Theory and Related Fields}
\bvolume{174}
\bpages{945--979}.
\end{barticle}
\endbibitem

\bibitem[\protect\citeauthoryear{Jing, Shao and Wang}{2003}]{JSW03}
\begin{barticle}[author]
\bauthor{\bsnm{Jing},~\bfnm{Bing-Yi}\binits{B.-Y.}},
  \bauthor{\bsnm{Shao},~\bfnm{Qi-Man}\binits{Q.-M.}} \AND
  \bauthor{\bsnm{Wang},~\bfnm{Qiying}\binits{Q.}}
(\byear{2003}).
\btitle{Self-normalized Cram{\'e}r-type large deviations for independent random
  variables}.
\bjournal{The Annals of probability}
\bvolume{31}
\bpages{2167--2215}.
\end{barticle}
\endbibitem

\bibitem[\protect\citeauthoryear{Jing, Wang and Zhou}{2015}]{jing2015cramer}
\begin{barticle}[author]
\bauthor{\bsnm{Jing},~\bfnm{Bing-Yi}\binits{B.-Y.}},
  \bauthor{\bsnm{Wang},~\bfnm{Qiying}\binits{Q.}} \AND
  \bauthor{\bsnm{Zhou},~\bfnm{Wang}\binits{W.}}
(\byear{2015}).
\btitle{Cram{\'e}r-type moderate deviation for studentized compound Poisson
  sum}.
\bjournal{Journal of Theoretical Probability}
\bvolume{28}
\bpages{1556--1570}.
\end{barticle}
\endbibitem

\bibitem[\protect\citeauthoryear{Krylov and Priola}{2010}]{KP08}
\begin{barticle}[author]
\bauthor{\bsnm{Krylov},~\bfnm{Nicolai~V.}\binits{N.~V.}} \AND
  \bauthor{\bsnm{Priola},~\bfnm{Enrico}\binits{E.}}
(\byear{2010}).
\btitle{Elliptic and parabolic second-order PDEs with growing coefficients}.
\bjournal{Communications in Partial Differential Equations}
\bvolume{35}
\bpages{1--22}.
\end{barticle}
\endbibitem

\bibitem[\protect\citeauthoryear{Lin and Bai}{2011}]{lin2011probability}
\begin{bbook}[author]
\bauthor{\bsnm{Lin},~\bfnm{Zhengyan}\binits{Z.}} \AND
  \bauthor{\bsnm{Bai},~\bfnm{Zhidong}\binits{Z.}}
(\byear{2011}).
\btitle{Probability inequalities}.
\bpublisher{Springer Science \& Business Media}.
\end{bbook}
\endbibitem

\bibitem[\protect\citeauthoryear{Majka, Mijatovi{\'c} and
  Szpruch}{2020}]{Majka2020}
\begin{barticle}[author]
\bauthor{\bsnm{Majka},~\bfnm{Mateusz~B.}\binits{M.~B.}},
  \bauthor{\bsnm{Mijatovi{\'c}},~\bfnm{Aleksandar}\binits{A.}} \AND
  \bauthor{\bsnm{Szpruch},~\bfnm{{\L}ukasz}\binits{{\L}.}}
(\byear{2020}).
\btitle{Nonasymptotic bounds for sampling algorithms without log-concavity}.
\bjournal{Annals of Applied Probability}
\bvolume{30}
\bpages{1534--1581}.
\end{barticle}
\endbibitem

\bibitem[\protect\citeauthoryear{Mao}{2008}]{Mao97}
\begin{bbook}[author]
\bauthor{\bsnm{Mao},~\bfnm{Xuerong}\binits{X.}}
(\byear{2008}).
\btitle{Stochastic differential equations and applications}.
\bpublisher{Horwood Publishing Limited}.
\end{bbook}
\endbibitem

\bibitem[\protect\citeauthoryear{Massart}{2007}]{Massart2003}
\begin{bbook}[author]
\bauthor{\bsnm{Massart},~\bfnm{Pascal}\binits{P.}}
(\byear{2007}).
\btitle{Concentration inequalities and model selection}
\bvolume{6}.
\bpublisher{Springer}.
\end{bbook}
\endbibitem

\bibitem[\protect\citeauthoryear{McLeish}{1974}]{McL1974}
\begin{barticle}[author]
\bauthor{\bsnm{McLeish},~\bfnm{Donald~L.}\binits{D.~L.}}
(\byear{1974}).
\btitle{Dependent central limit theorems and invariance principles}.
\bjournal{The Annals of Probability}
\bvolume{2}
\bpages{620--628}.
\end{barticle}
\endbibitem

\bibitem[\protect\citeauthoryear{Meyn and Tweedie}{1993}]{MeTw93}
\begin{barticle}[author]
\bauthor{\bsnm{Meyn},~\bfnm{Sean~P.}\binits{S.~P.}} \AND
  \bauthor{\bsnm{Tweedie},~\bfnm{Richard~L.}\binits{R.~L.}}
(\byear{1993}).
\btitle{Stability of Markovian processes III: Foster-Lyapunov criteria for
  continuous-time processes}.
\bjournal{Advances in Applied Probability}
\bvolume{25}
\bpages{518--548}.
\end{barticle}
\endbibitem

\bibitem[\protect\citeauthoryear{Meyn and Tweedie}{2009}]{MeTw93b}
\begin{bbook}[author]
\bauthor{\bsnm{Meyn},~\bfnm{Sean~P.}\binits{S.~P.}} \AND
  \bauthor{\bsnm{Tweedie},~\bfnm{Richard~L.}\binits{R.~L.}}
(\byear{2009}).
\btitle{Markov chains and stochastic stability}.
\bpublisher{Cambridge University Press}.
\end{bbook}
\endbibitem

\bibitem[\protect\citeauthoryear{Nyquist}{2017}]{Nyq17}
\begin{barticle}[author]
\bauthor{\bsnm{Nyquist},~\bfnm{Pierre}\binits{P.}}
(\byear{2017}).
\btitle{Moderate deviation principles for importance sampling estimators of
  risk measures}.
\bjournal{Journal of Applied Probability}
\bvolume{54}
\bpages{490--506}.
\end{barticle}
\endbibitem

\bibitem[\protect\citeauthoryear{Partington}{2004}]{partington04}
\begin{bbook}[author]
\bauthor{\bsnm{Partington},~\bfnm{Jonathan~R}\binits{J.~R.}}
(\byear{2004}).
\btitle{Linear operators and linear systems: an analytical approach to control
  theory}
\bvolume{60}.
\bpublisher{Cambridge University Press}.
\end{bbook}
\endbibitem

\bibitem[\protect\citeauthoryear{Revuz and Yor}{2013}]{revuz2013continuous}
\begin{bbook}[author]
\bauthor{\bsnm{Revuz},~\bfnm{Daniel}\binits{D.}} \AND
  \bauthor{\bsnm{Yor},~\bfnm{Marc}\binits{M.}}
(\byear{2013}).
\btitle{Continuous martingales and Brownian motion}
\bvolume{293}.
\bpublisher{Springer Science \& Business Media}.
\end{bbook}
\endbibitem

\bibitem[\protect\citeauthoryear{Roberts and Rosenthal}{2004}]{Rob04}
\begin{barticle}[author]
\bauthor{\bsnm{Roberts},~\bfnm{Gareth~O.}\binits{G.~O.}} \AND
  \bauthor{\bsnm{Rosenthal},~\bfnm{Jeffrey~S.}\binits{J.~S.}}
(\byear{2004}).
\btitle{General state space Markov chains and MCMC algorithms}.
\bjournal{Probability surveys}
\bvolume{1}
\bpages{20--71}.
\end{barticle}
\endbibitem

\bibitem[\protect\citeauthoryear{Roberts and Tweedie}{1996}]{Tweedie1996}
\begin{barticle}[author]
\bauthor{\bsnm{Roberts},~\bfnm{Gareth~O.}\binits{G.~O.}} \AND
  \bauthor{\bsnm{Tweedie},~\bfnm{Richard~L.}\binits{R.~L.}}
(\byear{1996}).
\btitle{Exponential convergence of Langevin distributions and their discrete
  approximations}.
\bjournal{Bernoulli}
\bvolume{2}
\bpages{341--363}.
\end{barticle}
\endbibitem

\bibitem[\protect\citeauthoryear{Shao}{1999}]{shao1999cramer}
\begin{barticle}[author]
\bauthor{\bsnm{Shao},~\bfnm{Qi-Man}\binits{Q.-M.}}
(\byear{1999}).
\btitle{A Cram{\'e}r type large deviation result for Student's t-statistic}.
\bjournal{Journal of Theoretical Probability}
\bvolume{12}
\bpages{385--398}.
\end{barticle}
\endbibitem

\bibitem[\protect\citeauthoryear{Shao}{2018}]{Shao18}
\begin{barticle}[author]
\bauthor{\bsnm{Shao},~\bfnm{Jinghai}\binits{J.}}
(\byear{2018}).
\btitle{Invariant Measures and Euler--Maruyama's Approximations of
  State-Dependent Regime-Switching Diffusions}.
\bjournal{SIAM Journal on Control and Optimization}
\bvolume{56}
\bpages{3215--3238}.
\end{barticle}
\endbibitem

\bibitem[\protect\citeauthoryear{Shao, Zhang and Zhang}{2018}]{SZZ20+}
\begin{barticle}[author]
\bauthor{\bsnm{Shao},~\bfnm{Qi-Man}\binits{Q.-M.}},
  \bauthor{\bsnm{Zhang},~\bfnm{Mengchen}\binits{M.}} \AND
  \bauthor{\bsnm{Zhang},~\bfnm{Zhuo-Song}\binits{Z.-S.}}
(\byear{2018}).
\btitle{Cram{\'e}r-type Moderate Deviation Theorems for Nonnormal
  Approximation}.
\bjournal{arXiv preprint arXiv:1809.07966}.
\end{barticle}
\endbibitem

\bibitem[\protect\citeauthoryear{Shao and Zhou}{2016}]{shao2016cramer}
\begin{barticle}[author]
\bauthor{\bsnm{Shao},~\bfnm{Qi-Man}\binits{Q.-M.}} \AND
  \bauthor{\bsnm{Zhou},~\bfnm{Wen-Xin}\binits{W.-X.}}
(\byear{2016}).
\btitle{Cram{\'e}r type moderate deviation theorems for self-normalized
  processes}.
\bjournal{Bernoulli}
\bvolume{22}
\bpages{2029--2079}.
\end{barticle}
\endbibitem

\bibitem[\protect\citeauthoryear{Teh, Thiery and Vollmer}{2016}]{Teh2016}
\begin{barticle}[author]
\bauthor{\bsnm{Teh},~\bfnm{Yee~Whye}\binits{Y.~W.}},
  \bauthor{\bsnm{Thiery},~\bfnm{Alexandre~H.}\binits{A.~H.}} \AND
  \bauthor{\bsnm{Vollmer},~\bfnm{Sebastian~J.}\binits{S.~J.}}
(\byear{2016}).
\btitle{Consistency and fluctuations for stochastic gradient Langevin
  dynamics}.
\bjournal{The Journal of Machine Learning Research}
\bvolume{17}
\bpages{193--225}.
\end{barticle}
\endbibitem

\bibitem[\protect\citeauthoryear{Tierney}{1994}]{Ti94}
\begin{barticle}[author]
\bauthor{\bsnm{Tierney},~\bfnm{Luke}\binits{L.}}
(\byear{1994}).
\btitle{Markov chains for exploring posterior distributions}.
\bjournal{the Annals of Statistics}
\bvolume{22}
\bpages{1701--1762}.
\end{barticle}
\endbibitem

\end{thebibliography}

\end{document}